\documentclass{amsart}

\numberwithin{equation}{section}

\usepackage{amssymb}
\usepackage{enumerate}
\usepackage{xspace}
\usepackage{hyperref}
\usepackage[latin1]{inputenc}
\usepackage[T1]{fontenc}
\usepackage{comment}

\newtheorem{thm}{Theorem}[section]
\newtheorem{lem}[thm]{Lemma}
\newtheorem{cor}[thm]{Corollary}

\newtheorem{prop}[thm]{Proposition}

\newtheorem{rem}[thm]{Remark}

\newcommand\B{{\mathcal B}}
\newcommand\Co{{\mathcal C}}
\newcommand\D{{\mathcal D}}

\newcommand\Lp{{\mathcal L}}

\newcommand\A{{\mathbb A}}
\newcommand\C{{\mathbb C}}

\newcommand\N{{\mathbb N}}
\newcommand\R{{\mathbb R}}
\newcommand\To{{\mathbb T}}

\newcommand\ve{\varepsilon}
\newcommand\vf{\varphi}

\newcommand\Id{\operatorname{{\bf Id}}}

\renewcommand{\div}{\operatorname{div}}
\newcommand{\tq}{\ |\ }

\newcommand{\norm}[1]{\left\| #1 \right\|}
\DeclareMathOperator{\Vol}{Vol}

\newcommand{\dd}{\, {\rm d}}
\newcommand{\ovB}{\overline{B}}
\renewcommand{\tilde}{\widetilde}
\newcommand{\spectr}{\ensuremath{\text{sp}}} 
\newcommand{\dist}{\operatorname{dist}}
\newcommand{\rank}{\operatorname{rank}}

\newcommand{\oo}{\overline{\omega}}
\newcommand{\rad}{\varrho}
\newcommand{\parti}{\vartheta}
\newcommand{\vectfield}{\mathcal{V}}
\newcommand{\chart}{\theta}
\DeclareMathOperator{\Jac}{Jac}
\newcommand{\smp}{t}

\begin{document}

\title{Banach spaces adapted to Anosov systems}

\author{S\'ebastien Gou\"ezel and Carlangelo Liverani}
\address{S\'ebastien Gou\"ezel\\
D\'epartement de math\'ematiques et applications\\
École Normale Sup\'erieure\\
45 rue d'Ulm, 75005 Paris, France.}
\email{sebastien.gouezel@ens.fr}
\address{Carlangelo Liverani\\
Dipartimento di Matematica\\
II Università di Roma (Tor Vergata)\\
Via della Ricerca Scientifica, 00133 Roma, Italy.}
\email{liverani@mat.uniroma2.it}
\thanks{We wish to thank V.Baladi, G.Keller and D.Ruelle for helpful
discussions. Last, but most of all, it is a pleasure to
acknowledge our debit to D.Dolgopyat. In fact, a couple of years
ago, while in Rome, he told the second named author that one
should ``project along the unstable direction'' and a few other
words. It is a fact that the above sentence is the key idea of the
present work. It is only to our detriment that it took us so long
to figure out the meaning and to understand how to implement it.
C.Liverani acknowledges the support of M.I.U.R. and the
hospitality of Paris seven and Georgia Institute of Technology
where part of the paper was written.}
\date{February 21, 2005}
\begin{abstract}
We study the spectral properties of the Ruelle-Perron-Frobenius
operator associated to an Anosov map on classes of functions with
high smoothness. To this end we construct anisotropic Banach
spaces of distributions on which the transfer operator has a small
essential spectrum. In the $\Co^\infty$ case, the essential
spectral radius is arbitrarily small, which yields a description
of the correlations with arbitrary precision. Moreover, we obtain
sharp spectral stability results for deterministic and random
perturbations. In particular, we obtain differentiability results
for spectral data (which imply differentiability of the SRB
measure, the variance for the CLT, the rates of decay for smooth
observable, etc.).
\end{abstract}
\maketitle

\section{introduction} \label{sec:intro}

The study of the statistical properties of Anosov systems dates
back almost half a century (\cite{Anos}) and many approaches have
been developed to investigate various aspects of the field (the
most historically relevant one being based on the introduction of
Markov partitions \cite{AS, Sin1, Bow1, Rue0}). At the same time
the type of questions and the precision of the results have
progressed through the years. In the last years the emphasis has
been on strong stability properties with respect to various types
of perturbations \cite{Bal}, dynamical zeta functions and related
smoothness issue (see \cite{Ki, Rue1, Dolgo1}). In the present
paper we present a new approach, improving on a previous partial
and still unsatisfactory one \cite{BKL}, that allows to obtain
easily a manifold of results (many of which new) and we hope will
reveal an even larger field of applicability. Indeed, the ideas in
\cite{BKL} have already been applied with success to some
partially hyperbolic situations (flows) \cite{Li3} and we expect
them to be applicable to the study of dynamical zeta functions.

The basic idea is inspired by the work on piecewise expanding
maps, starting with \cite{LY, Kel1} and the many others that
contributed subsequently (see \cite{Bal} for a nice review on the
subject). That is to study directly the transfer operator (often
called the Ruelle-Perron-Frobenius operator) on appropriate
functional spaces.\footnote{This, to our knowledge, has been the
first Markov-partition--free approach to the study of the
statistical properties of systems with sensitive dependence on
initial conditions.} For the case of smooth expanding maps, the
Sobolev spaces $W^{n,1}$, or the Banach spaces $\Co^n$, turn out
to be proper spaces where the transfer operator acts as a
\emph{smoothing} operator, \cite{Rue-d, Li3}.\footnote{The choice
of $\Co^n$, which requires a bit more work (the analogous of the
argument at the end of subsection \ref{sec:LY0} here) is the
choice generalized by the spaces we introduce in this paper.} In
turn, this implies that, on such spaces, the operator is
quasi-compact with an essential spectral radius exponentially
decreasing in $n$. The existence of a spectral gap and all kind of
statistical properties (exponential decay of correlations, central
limit theorem, meromorphic zeta functions, etc.) readily follow.

Unfortunately, for Anosov systems it is not helpful to consider
spaces of smooth functions -- on such spaces the spectral radius of
the transfer operator is larger than one --, it is necessary to
consider spaces of distributions. This was recognized in
\cite{Rugh1, Rugh2, Rugh3, Fried1} limited to the analytic case, and
in \cite{Li1} (only implicitly) and systematically in \cite{BKL} for
the $\Co^{1+\alpha}$ case. Nevertheless, the latter setting had
still several shortcomings. First of all, the Banach space was
precisely patterned on the invariant distributions of the systems,
which implied that transfer operators -- even of close maps -- where
studied on different spaces. This was a serious obstacle to
obtaining sharp perturbation results. Secondly, since in general the
invariant distributions are only H\"older, it was not possible to
have a scale of Banach spaces on which to study the influence of the
smoothness of the map on the spectrum.

Both such shortcomings are overcome in the present approach. The
spaces we introduce (partially inspired by \cite{Li1}) are still
related to the map one wishes to
study, but in a much loser way so that the operators associated to
nearby maps can be studied on the same space. In addition, we have
a scale of spaces that can be used to investigate smoothness
related issues (typically the dependence of the essential spectrum
on the smoothness of the map). In particular, if the map is
$\Co^\infty$, we obtain a description of the correlations of
$\Co^\infty$ functions with an arbitrarily small error term.

In addition, the present norms allow easier estimates of the size
of perturbations. This provides a very direct way of obtaining
sharp perturbations results which substantially generalize
the existing ones, e.g. \cite{BKL, Pol1, Rue1, Rue2}. For example,
in the $\Co^\infty$ case all the simple eigenvalues and all the
eigenspaces depend $\Co^\infty$ on the map. The same holds for the
variance in the CLT for a smooth zero average observable.

A further remarkable feature of the present approach is that,
unlike all the previous ones, its implementation does not depend
directly on subtle regularity properties of the foliations and of
the holonomies. This makes possible to have a much simpler and
{\sl self contained} treatment of the statistical properties of
the system and may lead to interesting generalizations in the
partially hyperbolic setting.

The paper is organized as follows. In the second section, we
introduce Banach spaces $\B^{p,q}$, explain why the transfer
operator acting on $\B^{p,q}$ has a spectral gap and illustrate
the stability results: the main ingredients are a compactness
statement (Lemma \ref{lem:comp}), a Lasota-Yorke type inequality
(Lemma \ref{lem:LY}) and the estimates on perturbations Lemmas
\ref{lem:mapdist} and \ref{lem:estrandom}. In Sections 3 and 4, we
describe more precisely the spaces $\B^{p,q}$ and prove in
particular that they are spaces of distributions. In Sections 5
and 6, which are the main parts of this article, we prove
respectively the aforementioned compactness statement and
Lasota-Yorke type inequality. In Section 7, we show how this
framework implies very precise stability results on the spectrum,
for deterministic and random perturbations. Section
\ref{sec:smoothpert} contains an abstract perturbation result
generalizing the setting of \cite{KL}, along the direction
adumbrated in \cite{Li3}, to cases where a control on the
smoothness is available. Section \ref{sec:smooth} shows that
smooth deterministic perturbations fit in the setting developed in
Section \ref{sec:smoothpert}. Finally, an appendix contains the
proof of an intuitive, but technical, result.

\section{The Banach spaces and the results} \label{sec:space}

For $q\geq 0$, let $\lfloor q\rfloor$ be its integer part. We
denote by $\bar \Co^q$ the set of functions which are $\lfloor
q\rfloor$ times continuously differentiable and whose $\lfloor
q\rfloor$-th derivative is H\"older continuous of exponent
$q-\lfloor q\rfloor$, if $q$ is not an integer. To fix notation,
in this paper we choose, for each $q\in\R_+$, a norm on
$\bar\Co^q$ functions so that $|\vf_1 \vf_2|_{\Co^q} \leq
|\vf_1|_{\Co^q} |\vf_2|_{\Co^q}$. We will denote by $\Co^q$ the
closure in $\bar \Co^q$ of the set of $\Co^\infty$ functions. It
coincides with $\bar \Co^q$ if $q$ is an integer, but is strictly
included in it otherwise. In any case, it contains $\bar \Co^{q'}$
for all $q'>q$.

Let $X$ be a $d$ dimensional $\Co^\infty$ compact connected
Riemannian manifold and consider an Anosov map
$T\in\Co^{r+1}(X,X)$ (for some real $r>1$). 
Write $d_s$ and $d_u$
for the stable and unstable dimensions. Let $\lambda>1$ be less
than the minimal expansion along the unstable directions, $\nu<1$
greater than the minimal contraction along the stable directions.
We will express the spectral properties of $T$ using the constants
$\lambda$ and $\nu$.

In Section \ref{sec:defleaves}, we will define a set $\Sigma$ of
admissible leaves. The elements of $\Sigma$ are small
$\Co^{r+1}$
embedded compact manifolds with boundary, of dimension $d_s$,
close to local stable manifolds.\footnote{The precise definition
of the set $\Sigma$ is given by \eqref{eq:def_Omega}.}

In what follows, if $v$ is a smooth vector field on an open subset
of $X$ and $f$ is a smooth function, then $vf$ will denote the
derivative of $f$ in the direction $v$. If $v_1,\dots,v_p$ are
smooth vector fields, then $v_1\dots v_pf$ will denote
$v_1(v_2(\dots (v_pf))\dots)$. We will sometimes write
$\prod_{i=1}^p v_i f$ for this expression, although it may be a
little misleading since the vector fields $v_i$ do not necessarily
commute.

We are now ready to introduce the relevant norms. When $W\in
\Sigma$, we will denote by $\Co_0^q(W,\R)$ the set of functions
from $W$ to $\R$ which belong to $\Co^q$ and vanish on a
neighborhood of the boundary of $W$, and by $\vectfield^{r}(W)$
the set of $\Co^{r}$ vector fields defined on a neighborhood of
$W$.

\relax For each $h\in\Co^r(X,\R)$ and $q\in \R_+$, $p\in\N$ with $p\leq
r$ (recall, $T$ is $\Co^{r+1}$ by definition), let\footnote{All
integrals are taken with respect to Lebesgue or Riemannian measure
except when another measure is explicitly mentioned.}
  \begin{equation}
  \label{eq:norm-st_ini}
  \|h\|^-_{p,q}:=
  \sup_{W \in \Sigma}\;\;\;
  \sup_{\substack{v_1,\dots,v_p \in \vectfield^{r}(W)\\ |v_i|_{\Co^{r}}
  \leq 1}}\;\;\;
  \sup_{\substack{\vf\in \Co_0^{q}(W,\R)\\
  |\vf|_{\Co^{q}}\leq 1}}\;\int_W  v_1\dots v_ph
  \cdot \vf.
  \end{equation}
It satisfies $\|h\|^-_{p,q'}\leq \|h\|^-_{p,q}$ if $q'\geq q$.
Define then the norms
  \begin{equation}\label{eq:norm-st}
  \|h\|_{p,q} = \sup_{0\leq k\leq p} \|h\|_{k,q+k}^- = \sup_{p'\leq
p, q'\geq q+p'} \|h\|^-_{p',q'}.
  \end{equation}
For example, if $X$ is the torus, the above norm is equivalent to
the one given by
  \begin{equation*}
  \sup_{|\alpha|\leq p}\;\;
  \sup_{W \in \Sigma}\;\;
  \sup_{\substack{\vf\in \Co_0^{q+|\alpha|}(W,\R)\\
  |\vf|_{\Co^{q+|\alpha|}}\leq 1}}\;\int_W  \partial^\alpha h
  \cdot \vf.
  \end{equation*}
Later on, in Section
\ref{sec:defleaves}, we will give an explicit description of the
norm \eqref{eq:norm-st} in coordinate charts. It will sometimes be
easier to work with the coordinate-free definition given in
\eqref{eq:norm-st} and sometimes with the explicit definition,
depending on what we are trying to prove.

It is easy to see that $\|\cdot \|_{p,q}$,
for $p\leq r$,
is a norm on $\Co^r(X,\R)$
(we will prove a more general result in Proposition
\ref{prop:Bpqdist}). Hence, we can consider the completion
$\B^{p,q}$ of $\Co^r(X,\R)$ with respect to this norm. Section
\ref{sect:descBpq} will be devoted to a description of this space.
We will see in particular that it is canonically a space of
distributions.

Since $\|h\|_{p-1,q+1} \leq \|h\|_{p,q}$, the embedding of
$\Co^r(X,\R)$ into $\B^{p,q}$ gives rise to a canonical map
$\B^{p,q} \to \B^{p-1,q+1}$, which is in fact compact:
\begin{lem}\label{lem:comp}
If $p+q<r$, the unit ball of $\B^{p,q}$ is relatively compact in
$\B^{p-1,q+1}$.
\end{lem}
The proof of Lemma \ref{lem:comp} is the content of Section
\ref{sec:comp}.

The rest of the paper consists in the
investigation of the properties of the transfer operator
$\Lp$ seen as an operator acting on the spaces $\B^{p,q}$.
As is well known, for each
$h\in\Co^{r}(X,\R)$, the transfer operator
$\Lp:\Co^{r}(X,\R)\to \Co^{r}(X,\R)$, defined by
duality by
\[
\int h\cdot u\circ T =: \int \Lp h \cdot u,
\]
is also given by
\[
\Lp h=(h |\det(DT)|^{-1})\circ T^{-1}.
\]
The key information on the action of $\Lp$ on $\B^{p,q}$ is
contained in the next lemma.
\begin{lem}\label{lem:LY}
For each $p\in\N$ and $q\geq 0$ satisfying $p+q< r$,
 $\Lp$ is a bounded operator on
$\B^{p,q}$.\footnote{That is, $\Lp$ can be extended to a bounded
operator on $\B^{p,q}$ that, with a mild abuse of notation, we
still call $\Lp$.} In addition, there exist $A_{p,q},B_{p,q}>0$
such that, for each $n\in\N$,
\begin{gather}
\|\Lp^n h\|_{0,q}\leq A_{0,q}\| h\|_{0,q}\,,\quad\text{for all }q<r; \label{eq1}
\\
\|\Lp^n h\|_{p,q}\leq A_{p,q}\max(\lambda^{-p},\nu^q)^{n}
\|h\|_{p,q}+B_{p,q}\|h\|_{p-1,q+1}\,,\quad\text{for all } p+q< r. \label{eq2}
\end{gather}
\end{lem}

The above Lemma is proven in section \ref{sec:LY}.

Lemmas \ref{lem:comp} and \ref{lem:LY} readily imply the basic result of the
paper:
\begin{thm}\label{thm:main}
If $p\in\N^*$ and $q\in \R_+^*$ satisfy $p+q< r$, then the
operator $\Lp:\B^{p,q}\to \B^{p,q}$ has spectral radius one. In
addition, $\Lp$ is quasicompact with essential spectrum
$\sigma_{\text{ess}}(\Lp) \subset\{z\in\C\;:\;|z|\leq
\max(\lambda^{-p},\nu^{q})\}$.

Moreover, the eigenfunctions corresponding to eigenvalues of modulus
$1$ are distributions of order $0$, i.e., measures.
If the map is topologically transitive, then one is a
simple eigenvalue, and no other eigenvalues of modulus one are
present.
\end{thm}
\begin{proof}
The first assertion follows from \eqref{eq2} since
$\|h\|_{p-1,q+1} \leq \|h\|_{p,q}$. The proof of the second is
completely standard and can be based, for example, on an argument
by Hennion after a spectral formula due to Nussbaum (see
\cite[Theorem 1]{BKL} for details).

The third is a consequence of the ergodic decomposition (see
\cite[Propositions 2.3.1 and 2.3.2]{BKL} for details).
\end{proof}

\begin{rem}\label{rem:notsogood}
If $\lambda^{-1}=\nu$ and $r$ is an even integer, then the optimal
choice of $p,q$ in Theorem \ref{thm:main}  is $p=r/2$ and
$q=r/2-\ve$ for some arbitrarily small $\ve>0$. Such a condition
is similar to Kitaev's requirement \cite{Ki} that, in our
language, reads $p=q=r/2$. Hence, our results are probably optimal
in this case. However, in the general case, our results are
limited by the fact that $p$ has to be an integer. It would
probably require significantly different and less elementary
techniques to allow $p\in \R_+$ as this
would require the definition of ``fractional derivatives''. See
\cite{Ba-l} for a very recent attempt, based on Fourier analysis, in a
very special case.
\end{rem}

In the $\Co^\infty$ case, Theorem \ref{thm:main} immediately
implies the following description of the correlations of
$\Co^\infty$ functions:
\begin{cor}
Assume that $T$ is $\Co^\infty$. Then there exist a sequence of
complex numbers
$\lambda_k$ such that $|\lambda_k|$ decreases to $0$, and
integers $r_k$ such that: for any $f,g: X\to \R$ of class $\Co^\infty$,
there exist numbers $a_k(f,g)$ with
  \begin{equation*}
  \int f \cdot g\circ T^n \sim \sum_{k=0}^\infty a_k(f,g) n^{r_k} \lambda_k^n
  \end{equation*}
in the following sense: for any $\ve>0$, let $K$ be such that
$|\lambda_K|<\ve$. Then
  \begin{equation*}
  \int f \cdot g\circ T^n = \sum_{k=0}^{K-1} a_k(f,g) n^{r_k}
  \lambda_k^n +o( \ve^n).
  \end{equation*}
\end{cor}

The second part of the paper focuses on the spectral stability for
a wide class of deterministic and random perturbations.

Let $U$ be a small enough neighborhood of $T$ in the $\Co^{r+1}$
topology. Consider a probability measure $\mu$ on a probability
space $\Omega$ and, for $\omega \in \Omega$, take $T_\omega \in U$
and $g(\omega, \cdot) \in  \Co^{p+q}(X, \R_+)$. Assume also that,
for all $x\in X$, $\int g(\omega,x) \dd\mu(\omega)=1$, and that
$\int |g(\omega, \cdot)|_{\Co^{p+q}(X,\R)} \dd\mu(\omega)
<\infty$. It is then possible to define a random walk in the
following way: starting from a point $x$, choose a diffeomorphism
$T_\omega$ randomly with respect to the measure $g( \omega,x)
\dd\mu(\omega)$, and go to $T_\omega(x)$. Then iterate this
process independently.

When $\Omega$ is a singleton and $g(\omega,x)=1$, then this is a
deterministic perturbation $T_\omega$ of $T$. Random perturbations
of the type discussed in \cite{BKL} can also be described in this
way.\footnote{To obtain the latter case set $\Omega=\To^d$, $\mu$ is
Lebesgue, $T_\omega(x)=Tx+\omega\mod 1$ and $g(\omega,x)=q_\ve(\omega,
Tx)$.}
Hence, this setting encompasses at the same time very general
deterministic and random perturbations of $T$.
Define the {\sl size} of the perturbation by
  \begin{equation}
  \label{eq:defDeltamug}
  \Delta(\mu,g):=\int |g(\omega,\cdot)|_{\Co^{p+q}(X,\R)}
  d_{\Co^{r+1}}(T_\omega,T) \dd\mu(\omega).
  \end{equation}
For definiteness, we will fix a large constant $A$ and assume
that, until the end of this paragraph, all the perturbations we
consider satisfy $\int |g(\omega,\cdot)|_{\Co^{p+q}(X,\R)} \leq
A$.

The transfer operator $\Lp_{\mu,g}$ associated to the previous
random walk is given by
  \begin{equation*}
  \Lp_{\mu,g} h (x)= \int_{\Omega} g(\omega, T_\omega^{-1}(x))
  \Lp_{T_\omega} h(x) \dd\mu(\omega)
  \end{equation*}
where $\Lp_{T_\omega}$ is the transfer operator associated to $T_\omega$.

In Lemma \ref{lem:mapdist} we show that $\Lp_T$ and
$\Lp_{\widetilde T}$ are $\|\cdot\|_{\B^{p,q}\to\B^{p-1,q+1}}$
close if $T$ and $\widetilde T$ are close in the $\Co^{r}$
topology. In turn, this implies that $\Lp_T$ and $\Lp_{\mu,g}$ are
close if $\Delta(\mu,g)$ is small, see \eqref{eq:kl_1}. In
addition, it is possible to show that the operators $\Lp_{\mu,g}$
satisfy a uniform Lasota-Yorke type inequality (Lemma
\ref{lem:estrandom}). These facts suffice to apply  \cite{KL} to
the present context, yielding immediately the strong perturbation
results described below greatly generalizing the results in
\cite{BKL}.

Fix any $\rad \in(\max(\lambda^{-p},\nu^{q}),1)$ and denote by
$\spectr(\Lp)$ the spectrum of $\Lp:\B^{p,q}\to\B^{p,q}$. Since
the essential spectral radius of $\Lp$ does not exceed
$\max(\lambda^{-p},\nu^{q})$, the set
$\spectr(\Lp)\cap\{z\in\C:|z|\geq \rad \}$ consists of a finite
number of eigenvalues $\lambda_1,\dots,\lambda_k$ of finite
multiplicity. Changing $\rad$ slightly we may assume that
$\spectr(\Lp)\cap\{z\in\C:\;|z|=\rad \}=\emptyset$. Hence there
exists $\delta_*<\rad-\max(\lambda^{-p},\nu^{q})$ such that
\begin{gather*}
  |\lambda_i-\lambda_j|>\delta_*\quad(i\neq j)\ ;\\
  \dist(\spectr(\Lp),\{|z|=\rad \})>\delta_*.
\end{gather*}

\begin{thm}
  \label{thm:spectral-stability}
  For each $\delta\in (0,\delta_*]$ and $\eta<1-\frac{\log \rad}
  {\log \max(\lambda^{-p},\nu^{q})}$,
  there exists
  $\ve_0$ such that for any perturbation $(\mu,g)$ of $T$
  satisfying $\Delta(\mu,g)\leq \ve_0$,
  \begin{enumerate}[a)]
  \item
  The spectral projectors
  \begin{equation}\label{eq:spectral-projectors}
  \begin{split}
    \Pi^{(j)}_{\mu,g}
    &:=
    \frac 1{2\pi\imath}\int_{\{|z-\lambda_j|
    =\delta\}}(z-\Lp_{\mu,g})^{-1}\,dz\,\\
    \Pi^{(\rad)}_{\mu,g}
    &:=
    \frac 1{2\pi\imath}\int_{\{|z|=\rad \}}(z-\Lp_{\mu,g})^{-1}\,dz\\
  \end{split}
  \end{equation}
  are well defined. We will denote by $\Pi^{(j)}_0$ and
  $\Pi^{(\rad)}_0$ the corresponding projectors for the unperturbed
  transfer operator $\Lp$.
  \item There is $K_1>0$ such that
    $\|\Pi^{(j)}_{\mu,g}-\Pi^{(j)}_0\|_{\B^{p,q}\to\B^{p-1,q+1}}\leq K_1\,
    \Delta(\mu,g)^{\eta}$ and
    $\|\Pi^{(\rad)}_{\mu,g}-\Pi^{(\rad)}_0\|_{\B^{p,q}\to\B^{p-1,q+1}}
    \leq K_1\,
    \Delta(\mu,g)^{\eta}$.
  \item $\rank(\Pi^{(j)}_{\mu,g})=\rank(\Pi^{(j)}_0)$.
  \item There is $K_2>0$ such that
    $\|\Lp_{\mu,g}^n\Pi^{(\rad)}_{\mu,g}\|_{p,q}\leq K_2\,{\rad}^n$
  for all $n\in\N$.
  \end{enumerate}
\end{thm}

If the perturbation enjoys stronger regularity properties, then
much sharper results can be obtained. Such results follow from a
generalization of \cite{KL}, along the lines of \cite{Li3}, that
can be found in Section \ref{sec:smoothpert}.

To keep the exposition simple let us restrict ourselves to
deterministic perturbations. Since $\Co^{r+1}(X,X)$ has naturally the
structure of a $\Co^{r+1}$ Banach manifold, it makes sense to consider
perturbations belonging to $\Co^s([-1,1],\Co^{r+1}(X,X))$, that is
curves $T_t$ of $\Co^{r+1}$ maps from $X$ to $X$ such that, when
viewed in coordinates, their firsts $s$ derivatives with respect to
$t$ are $\Co^{r+1}$ functions.

\begin{thm}\label{thm:smooth}
Let $T_t\in \Co^s([-1,1],\Co^{r+1}(X,X))$ and $T_0$ be an Anosov
diffeomorphism. Let $q>0$ and $p,s\in \N^*$ be such that $p+q+s<
r+1$. Then there exists $\delta_*>0$ such that, for all
$t\in[-\delta_*,\delta_*]$, the eigenvalues and eigenprojectors
$\lambda_i(t)$, $\Pi_i(t)$ associated to $\Lp_{T_t}$ with
$|\lambda_i(0)|> \max(\lambda^{-p},\nu^{q})$ satisfy:
\begin{enumerate}
\item if $\lambda_i(0)$ is simple, then $\lambda_i(t)\in\Co^{s-1}$;
\item $\Pi_i(t)\in \Co^{s-1}(\B^{p-1+s,q},\B^{p-1,q+s})$.
\end{enumerate}
\end{thm}
The above theorem is proven in Section \ref{sec:smooth} by showing
that the hypotheses of Theorem \ref{thm_1} hold in the present
context.

\begin{rem}
Notice that, in Theorem \ref{thm:smooth}, there is some limitation
to the differentiability, coming from the fact that $p$ is an
integer, namely $s <r$ (see also Remark \ref{rem:notsogood}). In
certain cases this can yield a weaker result than \cite{Pol1}
where, in the case $s=r+1$, it is proven that the eigenvalues are
$\Co^{r-1}$. Yet, \cite{Pol1} is limited to the peripheral
eigenvalues and gives much less information on the eigenspaces.
\end{rem}

\begin{rem}
Note that, although not explicitly stated, all the constants in
Theorems \ref{thm:spectral-stability} and \ref{thm:smooth} are
constructive and can actually be computed in specific examples (see
\cite{Li4} for a discussion of such issues).
\end{rem}

\begin{rem}
Beside the eigenvalues and the eigenprojectors, the above theory
implies results also for other physically relevant quantities. For
example, if $T_0$ is a transitive Anosov map and $f\in\Co^{r}$, let
$f_t=f-\int f\dd\mu_{SRB}(T_t)$.\footnote{By Theorem \ref{thm:main},
$\mu_{SRB}$ is simply the eigenvector of $\Lp_{T_t}$ associated to the
eigenvalue one.} It is well
known that $\frac 1{\sqrt n}\sum_{k=0}^{n-1} f_t\circ T_t^k$
converges in law to a Gaussian with zero mean and variance
\[
\sigma(t)^2=-\mu_{SRB}(f_t^2)+2\sum_{n=0}^{\infty}\mu_{SRB}(f_t\circ T_t^n
f_t)=-\mu_{SRB}(f_t^2)+2[(\Id-\widetilde\Lp_{T_t})^{-1}\mu_t](f_t)
\]
where $\mu_t(\vf):=\mu_{SRB}(f_t \vf)$ and $\widetilde \Lp_{T_t}$ is
the operator $\Lp_{T_t}$ restricted to the space
$V_0:=\{h\in\B^{p,q}\;:\; \int_Xh=0\}$. 
By the results of Sections
\ref{sec:smoothpert} and \ref{sec:smooth} it follows then that
$\sigma\in\Co^{s-1}$ and one can actually compute formulae for its Taylor
expansion up to order $s-1$.
\end{rem}

We conclude this section with a warning to the reader.
\begin{rem}
Through the paper we will use $C$ and $C_\alpha$ to designate  
generic constants
depending only on the map, the Banach spaces and, eventually, on
the parameter $\alpha$. Their actual numerical value can thus change from one
occurrence to the next.
\end{rem}

\section{Definition and properties of the admissible leaves}

\label{sec:defleaves}

Replacing the metric by an adapted metric \emph{à la Mather}
\cite{Mat}, we can assume that the expansion of $DT(x)$ along the
unstable directions is stronger than $\lambda$, the contraction
along the stable directions is stronger than $\nu$, and the angle
between the stable and unstable directions is everywhere arbitrarily
close to $\pi/2$. For small enough $\kappa$, we define the stable
cone at $x\in X$ by
  \begin{equation*}
  \Co(x)=\left\{ u+v \in T_x X \tq u\in E^s(x), v \perp E^s(x), \norm{v}
  \leq \kappa \norm{u} \right\}.
  \end{equation*}
If $\kappa$ is small enough,
$DT^{-1}(x)(\Co(x)\backslash \{0\})$ is included
in the interior of $\Co(T^{-1}x)$, and $DT^{-1}(x)$ expands the
vectors in $\Co(x)$ by $\nu^{-1}$.

There exists a finite number of $\Co^\infty$ coordinate charts
$\psi_1,\ldots,\psi_N$ such that $\psi_i$ is defined on a subset
$(-r_i,r_i)^d$ of $\R^d$ (with its standard euclidian norm), such
that
\begin{enumerate}
\item $D\psi_i(0)$ is an isometry.
\item $D\psi_i(0)\cdot \bigl(\R^{d_s}\times\{0\}\bigr)= E^s(\psi_i(0))$.
\item The $\Co^{r+1}$-norms of $\psi_i$ and its inverse are bounded by
$1+\kappa$.
\item There exists $c_i\in (\kappa,2\kappa)$
such that the cone $\Co_i=\{ u+v\in \R^d
\tq u\in \R^{d_s}\times \{0\}, v\in \{0\}\times \R^{d_u}, \norm{v}\leq
c_i \norm{u}\}$ satisfies the following property: for any $x\in
(-r_i,r_i)^d$, $D\psi_i(x) \Co_i \supset \Co(\psi_i x)$ and $DT^{-1}(
D\psi_i(x) \Co_i) \subset \Co(T^{-1}\circ \psi_i(x))$.
\item The manifold $X$ is covered by the open sets
$\bigl(\psi_i((-r_i/2,r_i/2)^d)\bigr)_{i=1\ldots N}$.
\end{enumerate}
It is easy to construct such a chart around any point of $X$, hence a finite
number of them is sufficient to cover the whole manifold by compactness.

Let $G_i(K)$ be the set of graphs of functions $\chi$ defined on a
subset of $(-r_i,r_i)^{d_s}$ and taking values in
$(-r_i,r_i)^{d_u}$, belonging to $\Co^{r+1}$, with $|D\chi| \leq
c_i$ (i.e., the tangent space to the graph of $\chi$ belongs to
the cone $\Co_i$) and with $|\chi|_{\Co^{r+1}} \leq K$.\footnote{A
function defined on an arbitrary subset $A$ of $\R^d$ is of class
$\Co^{r}$ if there exists a $\Co^{r}$ extension to an open
neighborhood of $A$. Its norm is the infimum of the norms of such
extensions.}

The following is a classical consequence of the uniform hyperbolicity
of $T$.
\begin{lem}
\label{lem:decnorm} If $K$ is large enough, then there exists $K'<K$
such that, for any $W\in G_i(K)$ and for any $1\leq j\leq N$, the
set $\psi_j^{-1}\circ T^{-1} \circ\psi_i(W)$ belongs to $G_j(K')$.
\end{lem}

If $\kappa$ is small enough, then
$\nu^{-1}>(1+\kappa)^2 \sqrt{1+4\kappa^2}$. Hence, there exists $A>0$ such
that
  \begin{equation}
  \label{eq:def_A}
  \frac{\nu^{-1}}{(1+\kappa)^2\sqrt{1+4\kappa^2}} (A-1)>A.
  \end{equation}
Take $\delta>0$ small enough so that $A\delta < \min(r_i)/6$.

We define an \emph{admissible graph} as a map $\chi$ defined on some
ball $\ovB(x,A\delta)$ included in $(-2r_i/3,2r_i/3)^{d_s}$, taking
its values in $(-2r_i/3,2r_i/3)^{d_u}$, with
$\operatorname{range}(\Id,\chi)\in G_i(K)$. Denote by $\Xi_i$ the
set of admissible graphs on $(-2r_i/3,2r_i/3)^{d_s}$.

Given an admissible graph $\chi\in\Xi_i$, we will call $\tilde
W:=\psi_i\circ(\Id,\chi)(\ovB(x,A\delta))$ the associated
\emph{full admissible leaf}
and $W:=\psi_i\circ(\Id,\chi)(\ovB(x,\delta))$ the \emph{admissible
leaf}.\footnote{Note that one can
talk about an admissible leaf only if it is given by
$\chi\in\Xi_i$, that is if there exists an associated full admissible
leaf.}

Let
  \begin{equation}
  \label{eq:def_Omega}
  \Sigma= \bigl\{ \psi_i \circ (\Id,\chi)(\ovB(x,\delta)) \tq
  \chi:\ovB(x,A\delta) \to \R^{d_u} \text{ belongs to } \Xi_i
  \bigr\}.
  \end{equation}
This is the set of admissible leaves.

We can use these admissible leaves to give another expression of
the norm \eqref{eq:norm-st} in coordinates. Set
  \begin{equation}\label{eq:norm-st1}
  \|h\|_{p,q}^{\sim} = \sup_{\substack{|\alpha|=p\\1\leq i\leq N}}\;
  \sup_{\substack{\chi :\ovB(x,A\delta) \to \R^{d_u}\\
  \chi \in \Xi_i }}\;
  \sup_{ \substack{ \vf \in \Co_0^{q}(\ovB(x,\delta),\R)
  \\ |\vf|_{\Co^{q}} \leq 1}}\; \int_{B(x,\delta)}
  [\partial^\alpha( h\circ \psi_i)] \circ (\Id,\chi) \cdot \vf.
  \end{equation}
and
  \begin{equation}
  \label{eq:norm-st2}
  \|h\|_{p,q}'=\sup_{0\leq k \leq p} \|h\|_{k,q+k}^{\sim}=
  \sup_{p'\leq p, q' \geq q+p'} \|h\|_{p',q'}^{\sim}.
  \end{equation}
The following lemma proves that this norm is equivalent to the
norm \eqref{eq:norm-st}, and gives a little bit more that will be
useful later.
\begin{lem}
\label{lem:equiv_norms}
If $p+q<r$, the norms $\|h\|_{p,q}$ and $\|h\|_{p,q}'$ are
equivalent. Moreover, there exists $C>0$ such that, for all $0\leq
k \leq p$, for all $f\in \Co^{q+k}$,
  \begin{equation}
  \label{eq:def_sharp}
  \sup_{W \in \Sigma}\;
  \sup_{\substack{v_1,\dots,v_{k} \in \vectfield^{q+k}(W)\\
  |v_i|_{\Co^{q+k}} \leq 1}} \;
  \sup_{\substack{\vf\in \Co_0^{q+k}(W,\R)\\
  |\vf|_{\Co^{q+k}}\leq 1}}\; \int_W  v_1\dots v_{k}(fh)
  \cdot \vf \leq C \|h\|'_{p,q} |f|_{\Co^{q+k}}.
  \end{equation}
\end{lem}
\begin{proof}
The inequality $\|h\|_{k,q+k}^{\sim} \leq C \|h\|_{k,q+k}^-$ is
trivial, since the images in the manifold of the coordinate vector
fields have a bounded $\Co^r$ norm. Hence, it is sufficient to
prove \eqref{eq:def_sharp} since, for $f=1$, it will imply
$\|h\|_{k,q+k}^- \leq C \|h\|'_{p,q}$.

We can without loss of generality work in one of the coordinate
charts $\psi_i$. Let $\chi$ be an admissible graph, and
$v_1,\dots, v_k$ be $\Co^{q+k}$ vector fields on a neighborhood of
the graph of $\chi$. Decomposing $v_j$ along the coordinate vector
fields, we can assume that $v_j= f_j
\partial_{\alpha(j)}$ where $f_j\in \Co^{q+k}$,
$|f_j|_{\Co^{q+k}}\leq 1$ and
$\alpha(j) \in \{1,\dots,d\}$. Hence,
  \begin{equation*}
  v_1 \dots v_k(f h)= \sum_{J,J_0,J_1,\dots,J_p} \Bigl( \prod_{j\in J}
  \partial_{\alpha(j)} h \Bigr) \Bigl( \prod_{j\in J_0}
  \partial_{\alpha(j)} f \Bigr)
  \Bigl( \prod_{j\in J_1} \partial_{\alpha(j)} f_1 \Bigr)
  \dots \Bigl( \prod_{j\in J_k} \partial_{\alpha(j)} f_k
  \Bigr),
  \end{equation*}
where we sum over all partitions $J,J_0,J_1,\dots,J_k$ of
$\{1,\dots,k\}$ such that $J_j \subset \{1,\dots,j-1\}$ for $j\geq
1$. Each term
  \begin{equation*}
  \int \biggl[ \Bigl( \prod_{j\in J}
  \partial_{\alpha(j)} h \Bigr)
  \Bigl( \prod_{j\in J_0} \partial_{\alpha(j)} f \Bigr)
  \Bigl( \prod_{j\in J_1} \partial_{\alpha(j)}
  f_1 \Bigr) \dots \Bigl( \prod_{j\in J_k} \partial_{\alpha(j)}
  f_k \Bigr) \biggr]\circ (\Id,\chi) \cdot \vf
  \end{equation*}
is an integral of $p'=|J|$ derivatives of $h$ against a test
function of differentiability $\Co^{q'}$ where $q'=\min_{0\leq
s\leq k}( q+k-|J_s|)\geq q+|J|$. Since $p'\leq p$ and $q'\geq
q+p'$, it is bounded by $C \|h\|_{p,q}'|f|_{\Co^{q+k}}$ by
\eqref{eq:norm-st2}.
\end{proof}
In the following, we will work indifferently with one expression
of the norm or the other and we will suppress the ``prime'' unless this
creates confusion.

The reason for integrating in
\eqref{eq:norm-st1}
only on admissible leaves, rather than on full admissible leaves,
is that the preimage of an admissible
leaf can be covered by a finite number of admissible leaves. We will
in fact need a slightly more precise result, conveniently expressed
in terms of the following notion. For $\gamma>1$, a
\emph{$\gamma$-admissible graph} is a map defined on a ball
$\ovB(x,\gamma A\delta)\subset
(-\frac{2r_i}{3\gamma},\frac{2r_i}{3\gamma})^{d_s}$, taking its
values in $(-\frac{2r_i}{3\gamma}, \frac{2r_i}{3\gamma})^{d_u}$,
whose graph belongs to $G_i(K)$. The corresponding
\emph{$\gamma$-admissible leaf} is $\psi_i \circ (\Id,\chi)(\ovB(x,
\delta/\gamma))$.
\begin{lem}
\label{lem:exists_partition}
There exists $\gamma_0>1$ satisfying the following property: for
any full admissible leaf $\tilde{W}$ and $n\in\N^*$, for any
$1\leq \gamma\leq \gamma_0$, there exist $\gamma$-admissible
leaves $W_1,\dots, W_\ell$, whose number $\ell$ is bounded by a
constant depending only on $n$, such that
\begin{enumerate}
\item $T^{-n}(W) \subset \bigcup_{j=1}^\ell W_j$.
\item $T^{-n}(\tilde{W}) \supset \bigcup_{j=1}^\ell W_j$.
\item There exists a constant $C$ (independent of $W$ and $n$) such that a
point of $T^{-n}\tilde{W}$ is contained in at most $C$ sets $W_j$.
\item There exist functions $\rho_1,\ldots, \rho_\ell$ of class
$\Co^{r+1}$ and compactly supported on $W_j$ such that $\sum \rho_j=1$
on $T^{-n}(W)$, and $| \rho_j|_{\Co^{r+1}} \leq C$.
\end{enumerate}
\end{lem}
\begin{proof}
Let $\chi:\ovB(x,A\delta) \to (-2r_i/3, 2r_i/3)^{d_u}$ be an
admissible graph. Let $W=\psi_i \circ (\Id,\chi)(\ovB(x,\delta))$
be the admissible leaf corresponding to $\chi$, and $\tilde W$ the
corresponding full admissible leaf.

Take $y\in \ovB(x,\delta)$ (so that $\ovB(y,(A-1)\delta) \subset
\ovB(x,A\delta)$) and $j$ such that $T^{-n} \circ \psi_i( y,\chi(y))
\in \psi_j \bigl( (-r_j/2,r_j/2)^d \bigr)$. Let $\pi: \R^d \to
\R^{d_s}$ be the projection on the first components. The map
$T^{-n}$ expands the distances by at least $\nu^{-n}$ along $\tilde
W$. The maps $\psi_i^{-1}$ and $\psi_j$ are $(1+\kappa)$-Lipschitz
and $|\pi(v)| \geq \frac{1}{\sqrt{1+4\kappa^2}} |v|$ when the vector
$v$ points in a stable cone $\Co_j$. Hence, the map
  $F:=\pi \circ \psi_j^{-1} \circ T^{-n} \circ \psi_i\circ (\Id, \chi)$
expands the distances by at least
$\frac{\nu^{-n}}{(1+\kappa)^2\sqrt{1+4\kappa^2} }$. If $\gamma$ is
close enough to $1$, then \eqref{eq:def_A} implies that the image by
$F$ of the ball $\ovB(y,(A-1)\delta)$ contains the ball $\ovB(F(y),
\gamma A \delta)$. We can then define a map $\chi_{F(y)} :
\ovB(F(y), \gamma A \delta) \to (-2r_j/(3\gamma),
2r_j/(3\gamma))^{d_u}$ such that its graph is contained in
$\psi_j^{-1}(T^{-n}\tilde W)$. In particular, $\chi_{F(y)}$ is a
$\gamma$-admissible graph, by Lemma \ref{lem:decnorm}.

We have shown that $T^{-n}W$ can be covered by $\gamma$-admissible
leaves. The lemma is then a consequence of \cite[Theorem
1.4.10]{Hor}.
\end{proof}
\begin{rem}
We will mostly use this lemma with $\gamma=1$, to get a covering of
$T^{-n}W$ by admissible leaves. However, in the study of
perturbations of $T$, we will need to use some $\gamma>1$.
\end{rem}

\section{Description of the space $\B^{p,q}$}

\label{sect:descBpq}

Take a covering of $X$ by sets of diameter at most $\delta$ and a
partition of unity subordinated to this covering. Using admissible
leaves supported in each of these sets, we easily check that there
exists a constant $C$ such that, for all $h\in \Co^r(X,\R)$ and
for all $\vf \in \Co^{q}(X,\R)$,
  \begin{equation*}
  \left| \int_X h\cdot \vf\right| \leq C \| h\|_{p,q}
  |\vf|_{\Co^{q}}.
  \end{equation*}
Passing to the completion, we obtain that any $h \in \B^{p,q}$
gives a distribution on $X$ of order at most $q$.\footnote{Here,
we are using the fact that there is a canonical measure on $X$,
namely the Riemannian measure. Otherwise, we would have to
distinguish between generalized functions and generalized
densities.} Denote by $\D'_{q}$ the set of distributions of order
at most $q$ with its canonical norm.
\begin{prop}
\label{prop:Bpqdist}
The map $\B^{p,q} \to \D'_{q}$ is a continuous injection.
\end{prop}
\begin{proof}
The continuity is trivial from the previous remarks.

Take $h\in \Co^r(X,\R)$. Let $\chi:\ovB(x,A\delta) \to
(-2r_i/3,2r_i/3)^{d_u}$ be an admissible graph and $|\alpha|\leq
p$. We can define a distribution $D_{\alpha,\chi}(h)$ of order
$q+|\alpha|$, on the ball $B(0,\delta)$, setting $\langle
D_{\alpha,\chi}(h),\vf \rangle = \int_{B(0,\delta)}
\partial^\alpha(h\circ \psi_i)(x+\eta,\chi(x+\eta)) \cdot
\vf(\eta) \dd\eta$. The map $h\mapsto D_{\alpha,\chi}(h)$ is
continuous for the $\|\cdot \|_{p,q}$-norm, whence it can be
extended to the space $\B^{p,q}$. The norm of an element $h$ of
$\B^{p,q}$ is by definition equal to the supremum of the norms of
the corresponding distributions $D_{\alpha,\chi}(h)$.

Assume that $\partial^\alpha=\partial_{j} \partial^\beta$. Let
$\chi_\ve$ be the admissible graph obtained by translating the
graph of $\chi$ of $\ve$ in the direction $x_j$. For $h\in \Co^r$,
the map $\ve \mapsto D_{\alpha,\chi_\ve} (h)$ is continuous. By
density, it is continuous for any $h\in \B^{p,q}$. Moreover,
  \begin{equation*}
  D_{\beta,\chi_\ve}(h)- D_{\beta,\chi}(h)= \ve \int_0^1 D_{\alpha,
  \chi_{t\ve}}(h) \dd t.
  \end{equation*}
Since $\ve \mapsto D_{\alpha,\chi_\ve}(h)$ is continuous, we
obtain that, for any $h\in \B^{p,q}$,
  \begin{equation}
  \label{eq:derDbeta}
  D_{\alpha,\chi}(h)=\lim_{\ve \to 0}
  \frac{D_{\beta,\chi_\ve}(h)- D_{\beta,\chi}(h)}{\ve}.
  \end{equation}

Take $h\in \B^{p,q}$ different from $0$. Then there exists an
admissible graph $\chi$ such that $D_{\emptyset,\chi}(h)\not=0$:
otherwise, \eqref{eq:derDbeta} would imply that all the
distributions $D_{\alpha,\chi}(h)$ vanish, which means that $h=0$
in $\B^{p,q}$. Since $\Co^\infty$ is dense in
$\Co^q$,\footnote{This is the only point in the paper for which it is useful to
consider $\Co^q$ instead of $\bar \Co^q$.} there exists $\vf\in
\Co^\infty$ such that $\langle D_{\emptyset,\chi}(h) ,\vf\rangle
\not =0$. Then, for any $\chi'$ close enough to $\chi$, we still
have
 $\langle D_{\emptyset,\chi'}(h) ,\vf\rangle \not =0$ by continuity.
Hence, we can construct a $\Co^\infty$ function $\tilde{\vf}$
supported on a neighborhood of the graph of $\chi$ such that
$\langle h, \tilde{\vf} \rangle \not=0$. Therefore, the
distribution given by $h$ is nonzero.
\end{proof}

\begin{rem}
There exist canonically defined maps $\B^{p,q} \to \B^{p-1,q}$ and
$\B^{p,q} \to \B^{p,q'}$ for $q'>q$, obtained by extending continuously the
canonical embedding of $\Co^r$ functions. Proposition
\ref{prop:Bpqdist} implies in particular that these maps are
injective.
\end{rem}

\begin{rem}
When $h$ is $\Co^r$ and $p+q< r$, then $\norm{h}_{p,q} \leq C
|h|_{\Co^r}$.
Hence, the embedding of $\Co^r(X,\R)$ in
$\B^{p,q}$ is continuous. Since $\Co^\infty(X,\R)$ is dense in
$\Co^r(X,\R)$ for the $\Co^r$-norm, it implies that
$\Co^\infty(X,\R)$ is dense in  $\B^{p,q}$. Hence, we could also
have obtained $\B^{p,q}$ by completing $\Co^\infty(X,\R)$.
\end{rem}

It is interesting to give explicit examples of nontrivial elements of
$\B^{p,q}$:

\begin{prop}
Let $W$ be a $\Co^{p+1}$-submanifold of dimension $d_u$ everywhere
transverse to the cones $D\psi_i(\Co_i)$ (e.g. a piece of unstable
manifold) and $\mu$ a $\Co^p$-density on $W$, with compact
support. Then the  distribution $\ell(\vf):= \int_W \vf \dd\mu$
belongs to $\B^{p,q}$.
\end{prop}
\begin{proof}
Without loss of generality we can assume that the manifold $W$
belongs to one chart $((-r_i,r_i)^d,\psi_i)$. We will work only in
such a chart, and omit the coordinate change $\psi_i$. The
pullback in this chart of the Riemannian measure is of the form
$\gamma(\eta,\xi) \dd\eta \dd\xi$, where $(\eta,\xi)\in
\R^{d_s}\times \R^{d_u}$.

The manifold $W$ is given by the graph of a $\Co^{p+1}$ function
$\zeta: \R^{d_u} \to \R^{d_s}$. The density of $\mu$ is then given
by a $\Co^p$ function $f: \R^{d_u} \to \R$ with compact support.

Let $f_\ve \in \Co^{\infty}$ be such that $|f-f_\ve|_{\Co^p} \leq
\ve$. Take also $\zeta_\ve \in \Co^\infty$ with $|\zeta
-\zeta_\ve|_{\Co^{p+1}} \leq \ve$. Let $\parti:\R^{d_s} \to \R_+$ be
a $\Co^\infty$ function supported in $B(0,1)$ and with $\int
\parti=1$. Set $\parti_\ve(\eta)=\frac{1}{\ve^{d_s}} \parti(\eta/\ve)$:
it is supported in $B(0,\ve)$ and has integral $1$. Let finally
$h_\ve(\eta,\xi)= \parti_\ve(\eta-\zeta_\ve(\xi))
f_\ve(\xi)\gamma(\xi,\eta)^{-1}$: it is a $\Co^\infty$ function,
and the corresponding distribution in $\D'_{q}$ is given by
  \begin{equation*}
  \vf \mapsto \int h_\ve(\eta,\xi)\vf(\eta,\xi)
  \gamma(\eta,\xi)\dd\eta\dd\xi
  = \int \parti_\ve(\eta-\zeta_\ve(\xi))f_\ve(\xi)
 \vf(\eta,\xi)\dd\eta\dd\xi.
  \end{equation*}
When $\ve\to 0$, this distribution converges to $\ell$.
Hence, the result will be
proved if we show that $\{h_\ve\}$ is a Cauchy sequence in
$\B^{p,q}$.

Take $\alpha$ with $|\alpha|\leq p$. Then one has
  \begin{equation*}
  \partial^\alpha h_\ve(\eta,\xi)= \sum_{\beta\leq \alpha}
 (\partial^\beta \parti_\ve)(\eta-\zeta_\ve(\xi))
  g_{\alpha,\beta,\ve}(\xi)
  \end{equation*}
where the function $g_{\alpha,\beta,\ve}$ is in $\Co^{\infty}$ and
converges in $\Co^{p-|\alpha|+|\beta|}$ to a function
$g_{\alpha,\beta,0}$ when $\ve \to 0$.

Let $\chi$ be an admissible graph and $\vf$ a $ \Co^{q+|\alpha|}$
test function with $|\vf|_{\Co^{q+|\alpha|}}\leq 1$. Then
  \begin{equation}
  \label{eq:Cauchy}
  \int \vf(\eta) \partial^\alpha h_\ve (\eta, \chi(\eta))
  = \sum_\beta \int \vf(\eta) (\partial^\beta
  \parti_\ve)(\eta-\zeta_\ve(\chi(\eta))
  g_{\alpha,\beta,\ve}(\chi(\eta)).
  \end{equation}
Since $W$ is everywhere transverse to the cone $\Co_i$, the map
$\theta_\ve : \eta \mapsto \eta-\zeta_\ve(\chi(\eta))$ is a
$\Co^{p+1}$ diffeomorphism, and it converges when $\ve \to 0$ to
$\theta_0$. Using this change of coordinates in \eqref{eq:Cauchy},
integrating by parts and given the fact that $\parti_\ve$ is a
$\Co^\infty$ mollifier, we obtain that \eqref{eq:Cauchy} converges
when $\ve \to 0$. Moreover, the speed of convergence is independent
of the graph $\chi$ or the test function $\vf$, since all norms are
uniformly bounded. Hence, $\norm{h_\ve - h_{\ve'}}_{p,q} \to 0$ when
$\ve, \ve' \to 0$, i.e., $h_\ve$ is a Cauchy sequence in $\B^{p,q}$.
\end{proof}

\section{Compactness} \label{sec:comp}

This paragraph is devoted to the proof of Lemma \ref{lem:comp}.
We will work only in coordinate charts, using
in an essential way the linear structure to interpolate between
admissible leaves.

Fix $\ve>0$. Take $1\leq i\leq N$. Since $p+q<r$, the injection $\Co^{r+1} \to
\Co^{p+q}$ is compact. Therefore, there exists a finite number of
admissible graphs $\chi_1,\ldots,\chi_s$ defined on balls
$\ovB(x_1,A\delta),\ldots, \ovB(x_s,A\delta)$ such that any other
admissible graph $\chi$ defined on a ball $\ovB(x,A\delta)$ is at
a distance at most $\ve$ of some $\chi_j$, in the sense that
$|x-x_j|\leq \ve$ and $|\eta \mapsto \chi(x+\eta)
-\chi_j(x_j+\eta)|_{\Co^{p+q}(\ovB(0,\delta),\R^{d_u})} \leq \ve$.

Take $\alpha$ with $|\alpha|=p-1$ and $\vf\in
\Co^{p+q}_0(\ovB(0,\delta),\R)$ with $|\vf|_{\Co^{p+q}} \leq 1$.
Write
  \begin{equation*}
  f_t(\eta)=(x_j + \eta + t(x-x_j) , \chi_j(x_j+\eta) +t
  (\chi(x+\eta)-\chi_j(x_j+\eta))).
  \end{equation*}
Write also $F(z)=\partial^\alpha (h\circ \psi_i)(z)$.
Then, for $\eta \in B(0,\delta)$,
  \begin{multline*}
  \partial^\alpha (h\circ \psi_i) (x+\eta,\chi(x+\eta))
  -
  \partial^\alpha (h\circ \psi_i) (x_j+\eta, \chi_j(x_j+\eta))
  =F( f_1(\eta))-F(f_0(\eta))
  \\
  =\int_0^1 DF (f_t(\eta))\cdot (x-x_j, \chi(x+\eta)-\chi_j(x_j+\eta))\dd t.
  \end{multline*}
Hence,
  \begin{multline*}
  \int \partial^\alpha (h\circ \psi_i) (x+\eta,\chi(x+\eta))
  \vf(\eta) \dd\eta
  - \int \partial^\alpha (h\circ \psi_i) (x_j+\eta,\chi_j(x_j+\eta))
  \vf(\eta)\dd \eta
  \\
  = \int_0^1
  \left(\int DF (f_t(\eta)) (x-x_j, \chi(x+\eta)-\chi_j(x_j+\eta))
  \vf(\eta) \dd \eta \right) \dd t.
  \end{multline*}
When $t$ is fixed, the last integral is an integral along the
graph given by $f_t$. This graph is admissible since it is an
interpolation between two admissible graphs (here, the fact that
the cone $\Co_i$ is constant is essential). Since $|x-x_j| \leq
\ve$ and $| \chi(x+\eta)-\chi_j(x_j+\eta)|_{\Co^{p+q}} \leq \ve$, 
this term can be estimated by $C \ve
\norm{h}^\sim_{p,q+p}$. We have proved that
  \begin{multline*}
  \|h\|^\sim_{p-1,q+p}=
  \sup_{\substack{|\alpha|=p\\1\leq i\leq N}} \;\sup_{\substack{\chi
  :\ovB(x,A\delta) \to \R^{d_u}\\
  \chi \in \Xi_i }}\;
  \sup_{\substack{\vf \in \Co_0^{p+q}(\ovB(x,\delta),\R)\\
  |\vf|_{\Co^{p+q}}\leq 1}}
  \int_{B(x,\delta)} \partial^\alpha(h\circ \psi_i)\circ
  (\Id,\chi)\cdot \vf
  \\
  \leq C\ve \norm{h}^\sim_{p,q+p} + \sup_{|\alpha|=p}\; \sup_{1\leq k
  \leq s}\;
  \sup_{\substack{\vf \in \Co_0^{p+q}(\ovB(x_k,\delta),\R)\\
  |\vf|_{\Co^{p+q}}\leq 1}}\;\;
  \int_{B(x_k,\delta)} \partial^\alpha(h\circ \psi_i)\circ
  (\Id,\chi_k)\cdot \vf,
  \end{multline*}
i.e., we have only a finite number of admissible graphs to
consider.

In the following, we work with one graph $\chi=\chi_k$. The set of
functions $\vf \in \Co_0^{p+q}(\ovB(x_k,\delta),\R)$ with
$|\vf|_{\Co^{p+q}\leq 1}$ is relatively compact for the
$\Co^{p+q-1}$ topology. Hence, there exists a finite set of
functions $\vf_1,\ldots,\vf_k$ which are $\ve$-dense. For any
$\vf$ as above, there exists $j$ such that
$|\vf-\vf_j|_{\Co^{p+q-1}} \leq \ve$. Since  $\|
h\|_{p-1,q+p-1}^\sim \leq \|h\|_{p,q}$,
  \begin{equation*}
  \int \partial^\alpha(h\circ \psi_i)\circ (\Id, \chi)\cdot
  \vf
  \leq \int \partial^\alpha(h\circ \psi_i)\circ (\Id, \chi)\cdot \vf_j
  +C \ve \| h\|_{p,q}.
  \end{equation*}

To summarize: we have proved the existence of a finite number of continuous
linear forms $\nu_1,\ldots,\nu_\ell$ on $\B^{p,q}$ such that, for any
$h\in \B^{p,q}$,
  \begin{equation*}
  \| h\|^\sim_{p-1,q+p} \leq C \ve \|h\|_{p,q} + \sup | \nu_i(h)|.
  \end{equation*}
This immediately implies the compactness we are looking for.

\section{Lasota-Yorke type inequality}\label{sec:LY}
This section is devoted to the proof of Lemma \ref{lem:LY}.

\subsection{Proof of \eqref{eq1}}\label{sec:LY0}
By density it suffices to prove it for $h\in\Co^r$.

Take $W\in \Sigma$ and $n\in \N^*$. Let $\vf \in \Co^q_0(W,\R)$
satisfy $|\vf|_{\Co^q}\leq 1$. Let $\rho_1, \ldots, \rho_\ell$ be
the partition of unity on $T^{-n}W$ given by Lemma
\ref{lem:exists_partition} (for $\gamma=1$), and
$W_1,\ldots,W_\ell$ the corresponding admissible leaves. Let $h_n=
h\cdot |\det DT^n|^{-1}$, then
  \begin{equation*}
  \int_{W} \Lp^n h \cdot \vf
  = \int_{T^{-n}W} h_n\cdot \vf\circ T^n\cdot J_W T^n
  \end{equation*}
where $J_W T^n$ is the jacobian of $T^n: T^{-n}W \to W$. Using the
partition of unity,
  \begin{equation}
  \label{eq:intermed0}
  \int_{W} \Lp^n h \cdot \vf=
  \sum_{j=1}^\ell \int_{W_j}
  h_n\cdot \vf\circ T^n\cdot J_W T^n\cdot \rho_j.
  \end{equation}
The function $\vf_j:=\vf\circ T^n \cdot \rho_j$ is compactly
supported on the admissible leaf $W_j$, and belongs to $\Co^q$.
Using the definition of the $\|\cdot \|_{0,q}$ norm along $W_j$
yields
  \begin{equation}
  \label{eq:intermed2}
  \left|\int_{W_j}
  h_n \cdot \vf_j\cdot J_W T^n  \right|
  \leq C \norm{h}_{0,q}
  \bigl| |\det DT^n|^{-1}\cdot \vf_j\cdot J_W T^n\bigr|_{\Co^{q}(W_j)}.
  \end{equation}
We will use repeatedly the following distortion lemmas:
\begin{lem}
\label{lem:distortion}
Let $\tilde{W}$ be a full admissible leaf and $W'$ an admissible
leaf contained in $T^{-n}\tilde{W}$. Let $1\leq s\leq r$. Let
$g_0,\dots,g_{n-1}$ be strictly positive $\Co^s$ functions on
$W',\dots, T^{n-1}(W')$ and $L>0$ be such that, for any $x\in
T^i(W')$, the $\Co^s$ norm of $g_i$ is bounded on an neighborhood
of $x$ by $L g_i(x)$. Then
  \begin{equation*}
  \forall x \in W', \quad \left| \prod_{i=0}^{n-1} g_i \circ T^i
  \right|_{\Co^s(W')} \leq C e^{CL} \prod_{i=0}^{n-1} g_i \circ T^i(x)
  \end{equation*}
for some constant $C$
depending only on the map $T$.
\end{lem}
\begin{proof}
Using the assumption on the $\Co^s$ norm of $g_i$ and the uniform
contraction of  $T$ along $W'$, it is easy to check that $\left|
\prod_{i=0}^{n-1} g_i \circ T^i \right|_{\Co^s(W')} \leq C L^s
\left| \prod_{i=0}^{n-1} g_i \circ T^i\right|_{\Co^0(W')}$.

The differential of the function $\log\left( \prod_{i=0}^{n-1} g_i
\circ T^i \right)$ is also bounded by $CL$, whence, for any
$x,y\in W'$,
  \begin{equation*}
  \prod_{i=0}^{n-1} g_i \circ T^i(x) \leq
  C e^{CL} \prod_{i=0}^{n-1} g_i \circ T^i(y).
  \qedhere
  \end{equation*}
\end{proof}

\begin{lem}
\label{lem:sumdet} There exists $C>0$ such that, for each $n\in\N$,
holds true
  \[
  \sum_{j\leq \ell}\left| |\det DT^n|^{-1}\right|_{\Co^{r}(W_j)}
  \cdot \left|J_W T^n\right|_{\Co^{r}(W_j)}\leq C.
  \]
\end{lem}
\begin{proof}
Lemma \ref{lem:distortion} applies to estimate $\left|\det
DT^n|^{-1} \right|_{\Co^{r}(W_j)}$ and $|J_W T^n|_{\Co^r(W_j)}$.
For any $x\in W_j$,
  \begin{equation*}
  \left| |\det DT^n|^{-1}\right|_{\Co^{r}(W_j)}
  \cdot \left|J_W T^n\right|_{\Co^{r}(W_j)}
  \leq C  |\det DT^n|^{-1}(x) J_W T^n(x).
  \end{equation*}
In particular,
  \begin{multline*}
  \left| |\det DT^n|^{-1}\right|_{\Co^{r}(W_j)}
  \cdot \left|J_W T^n\right|_{\Co^{r}(W_j)}
  \\
  \leq C \int_{W_j}|\det DT^n|^{-1} J_W T^n
  = C \int_{T^n (W_j)} |\det(DT^{-n})|.
  \end{multline*}
By Lemma \ref{lem:exists_partition}, the sets $T^n(W_j)$ are
contained in $\tilde{W}$ and have a bounded number of overlaps.
Let us consider the thickening $Z:=\bigcup_{ x\in \tilde{W}}
W^u_\rho(x)$, where $W^u_\rho(x)$ is the ball of size $\rho$ in
the unstable manifold through $x$. By usual distortion estimates,
  \begin{align*}
  \sum_{j\leq\ell} \int_{T^n (W_j)} |\det(DT^{-n})|&\leq C\int_{\tilde{W}}
  |\det(DT^{-n})| \\
  &\leq C \rho^{-d_u} \int_{Z}|\det DT^{-n}| = C \Vol(T^{-n}Z)\leq
  C. \qedhere
  \end{align*}
\end{proof}
Since $T^n$ is uniformly contracting along $T^{-n}(W)$, we have
$\left|\vf \circ T^n\right|_{\Co^{q}(W_j)} \leq C |\vf|_{\Co^{q}}
\leq C$, and $\left|\rho_j\right|_{\Co^{q}(W_j)}$ is uniformly
bounded by Lemma \ref{lem:exists_partition}. This, together with
\eqref{eq:intermed0} and \eqref{eq:intermed2} and Lemma
\ref{lem:sumdet}, concludes the proof of \eqref{eq1}.

\subsection{Proof of \eqref{eq2}}
Let $p\in \N$ and $q\geq 0$ satisfy $p+q<r$.

\begin{lem}
\label{main_lemma_LY}
There exists a constant $C$ such that, for each $n\in \N$, there
exists $C_n>0$ with
  \begin{equation*}
  \forall 0\leq \smp <p,\
  \norm{\Lp^n h}_{\smp,q+\smp}^- \leq C (\nu^q)^n \norm{h}_{p,q}+C_n \norm{h}_{p-1,q+1}
  \end{equation*}
and
  \begin{equation*}
  \norm{\Lp^n h}_{p,q+p}^- \leq C \max(\lambda^{-p},\nu^{q})^n
  \norm{h}_{p,q} + C_n \norm{h}_{p-1,q+1}.
  \end{equation*}
\end{lem}
\begin{proof}
We prove the lemma by induction over $\smp$. So, take $0\leq \smp
\leq p$ and assume that the conclusion of the lemma holds for all
$\smp'<\smp$.

Let $W$ and $\tilde{W}$ be an admissible leaf and the
corresponding full admissible leaf. As before, we will use Lemma
\ref{lem:exists_partition} (with $\gamma=1$) to write
$T^{-n}W\subset \bigcup_j W_j$ and denote by $\rho_j$ the
corresponding partition of unity given by Lemma
\ref{lem:exists_partition}.

Let $v_1,\dots,v_\smp \in \vectfield^r(W)$ with $|v_i|_{\Co^r}\leq
1$, and $\vf\in \Co^{\smp+q}_0(W)$ with
$|\vf|_{\Co^{\smp+q}(W)}\leq 1$. Writing $h_n=h\cdot |\det DT^n|$
as above, we want to prove that
  \begin{equation}
  \label{eq:main_to_prove}
  \left|\int_W v_1\dots v_\smp (h_n \circ T^{-n}) \vf \right|
  \\
  \leq \begin{cases}
  C (\nu^q)^n \norm{h}_{p,q} + C_n
  \norm{h}_{p-1,q+1} &\text{if }\smp<p.
  \\
  C \max(\lambda^{-p},\nu^{q})^n \norm{h}_{p,q} + C_n
  \norm{h}_{p-1,q+1} & \text{if }\smp=p.
  \end{cases}
  \end{equation}
The main idea of the proof will be to decompose each $v_i$ as a
sum $v_i=w_i^u+w_i^s$ where $w_i^s$ is tangent to $W$, and $w_i^u$
``almost'' in the unstable direction. We will then get rid of
$w_i^s$ by an integration by parts, and use the fact that $w_i^u$
is contracted by $DT^{-n}$ to conclude.

Let $\Psi$ and  $\Psi_j$ be coordinates charts with uniformly
bounded $\Co^{r+1}$ norms such that the images in the charts of
$\tilde W$ and $\tilde W_j$ are contained in $\R^{d_s}\times
\{0\}$. Let $\bar v_i(y)=D\Psi(\Psi^{-1}y)v_i(\Psi^{-1}y)$, it has
still a bounded $\Co^r$ norm. The integral in
\eqref{eq:main_to_prove} can be written as
  \begin{equation*}
  \int_{\Psi(W)} \bar v_1 \dots \bar v_\smp
  \bigl(h_n \circ T^{-n}\circ \Psi^{-1}\bigr)
  \cdot \vf\circ \Psi^{-1}\cdot \Jac(\Psi^{-1}).
  \end{equation*}
We will work in this coordinate chart and, with a small abuse of
notations, omit $\Psi$ in the formulas. Since $\Jac(\Psi^{-1})$
has a bounded $\Co^r$ norm, we may also replace $\vf$ with
$\vf\circ \Psi^{-1} \cdot \Jac(\Psi^{-1})$. We will also work in
the charts $\Psi_j$, and omit them as well in the formulas.
\begin{rem}
Note that, in the coordinate charts $\Psi, \Psi_j$ the manifolds $W,
W_j$ are $\Co^\infty$. This will be used extensively in the following.
\end{rem}

Without loss of generality, we can assume that each $\bar v_i$ is
of the form $f_i \partial_{\alpha(i)}$ where $f_i$ is bounded in
$\Co^r$ and $\partial_{\alpha(i)}$ is one of the coordinate vector
fields. In $\bar v_1\dots \bar v_\smp(h_n \circ T^{-n})$, if we
differentiate at least one of the functions $f_i$, we obtain an
integral of $\smp'<\smp$ derivatives of $h_n \circ T^{-n}$ against
a function in $\Co^{q+\smp'}$. Hence, it is bounded by $C
\norm{\Lp^n h}^-_{\smp',q+\smp'}$, which has already been
estimated in the induction. For the remaining term (where no $f_i$
has been differentiated), we can replace $\vf$ by $f_1\dots
f_\smp\cdot \vf$ and assume that $\bar v_i =\partial_{\alpha(i)}$.
In particular, $\bar v_i$ is well defined and smoth on a neighborhood of $\tilde 
W$.

The decomposition of $\bar v_i$ as $w_i^u+w_i^s$ on a neighborhood of
$T^n(W_j)$ is given by the following technical lemma:

\begin{lem}
\label{distortion}
Let $v$ be a vector field on a neighborhood of $\tilde W$ with
$|v|_{\Co^{r+1}}\leq 1$. Then there exist $\Co^{r+1}$ vector fields
$w^u$ and $w^s$ on a neighborhood $U$ of $T^n(W_j)$, satisfying:
\begin{itemize}
\item for all $x\in T^n(W_j)$, $w^s(x)$ is tangent to $T^n(W_j)$.
\item $|w^s|_{\Co^{r+1}(U)} \leq C_n$ and
$|w^u|_{\Co^{r+1}(U)}\leq C_n$, where $C_n$ is a constant that
may depend on $n$.
\item $|w^s \circ T^n|_{\Co^{r}(W_j)} \leq C$.
\item $|DT^n(x)^{-1}w^u(T^n x)|_{\Co^{p+q}(T^{-n}U)} \leq C
\lambda^{-n}$.
\end{itemize}
\end{lem}

The idea of the lemma is to decompose the tangent space at
$y=T^n(x)\in T^n(W_j)$ as the sum of the tangent space to $W$, and
the image of the vertical direction $\{0\}\times \R^{d_u}$ under
$DT^n(x)$. The decomposition $v=w^u+w^s$ is then obtained by
projecting $v$ along these two directions. The estimates on
$|w^s\circ T^n|$ and $|DT^n(x)^{-1}w^u(T^n x)|$ are then
consequences of the smoothing properties of $T^n$ along $W_j$.
This naive idea works well when $p+q\leq r-1$, but it yields only
$C^r$ vector fields $w^u$ and $w^s$, which is not sufficient for
our purposes if $r-1<p+q<r$. Hence, the rigorous proof of Lemma
\ref{distortion} involves additional regularization steps. Since
it is purely technical, it will be deferred to Appendix
\ref{appendix_distortion}.

Take some index $j$, we will estimate
  \begin{equation}
  \label{eq:to_estimate}
  \int_{T^n(W_j)} \bar v_1 \dots \bar v_\smp (h_n \circ T^{-n})
  \cdot \vf \cdot \rho_j
  \circ T^{-n}.
  \end{equation}
As in Lemma \ref{distortion}, write $\bar v_i= w_i^u+w_i^s$. Then
\eqref{eq:to_estimate} is equal to
  \begin{equation*}
  \sum_{\sigma \in \{s,u\}^\smp}
  \int_{T^n(W_j)} w_1^{\sigma_1} \dots w_\smp^{\sigma_\smp} (h_n \circ T^{-n})
  \cdot \vf \cdot \rho_j
  \circ T^{-n}.
  \end{equation*}
Take $\sigma \in \{s,u\}^\smp$, and let $k=\#\{i \tq
\sigma_i=s\}$. Let $\pi$ be a permutation of $\{1,\dots,\smp\}$
such that $\pi\{1,\dots, k\}= \{i \tq \sigma_i=s\}$. Then
  \begin{multline*}
  \int_{T^n(W_j)} w_1^{\sigma_1} \dots w_\smp^{\sigma_\smp} (h_n\circ
  T^{-n}) \cdot \vf \cdot \rho_j
  \circ T^{-n}
  =\\
  \int_{T^n(W_j)} \prod_{i=1}^k w_{\pi(i)}^s \prod_{i=k+1}^\smp
  w_{\pi(i)}^u (h_n\circ T^{-n})  \cdot \vf \cdot \rho_j
  \circ T^{-n} + O(\|h\|_{p-1,q+1}).
  \end{multline*}
Namely, the commutator of two $\Co^{r+1}$ vector fields is a
$\Co^{r}$ vector field. Hence, if we exchange two vector fields,
the difference is bounded by $C_n\|h\|_{p-1,q+1}$.

We integrate by parts with respect to the vector fields
$w_{\pi(i)}^s$: they are tangent to the manifold $W$, whence
$\int_W vf \cdot g= -\int_W f \cdot vg + \int_W fg \cdot \div v$.
Since $w_{\pi(i)}^s$ is $\Co^{r+1}$ and the manifold $W$ is
$\Co^{\infty}$ with a $\Co^{\infty}$ volume form (here, we use the
fact that we work in a coordinate chart for which $W\subset
\R^{d_s}\times\{0\}$), the divergence terms are bounded by
$C_n\|h\|_{p-1,q+1}$. We get an integral $\int \prod_{k+1}^\smp
w_{\pi(i)}^u (h_n\circ T^{-n}) \cdot \prod_{k}^1 w_{\pi(i)}^s
\left( \vf \cdot \rho_j \circ T^{-n}\right)$. If we use one of the
vector fields $w_{\pi(i)}^s$ to differentiate $\rho_j \circ
T^{-n}$, then we obtain an integral of $\smp-k$ derivatives of $h$
against a function in $\Co^{q+1 +(\smp-k)}$, which is again bounded by
$C_n \|h\|_{p-1,q+1}$. Hence,
  \begin{multline}
  \label{eq:kok}
  \int_{T^n(W_j)} \prod_{i=1}^k w_{\pi(i)}^s \prod_{i=k+1}^\smp
  w_{\pi(i)}^u (h_n \circ T^{-n}) \cdot \vf \cdot \rho_j
  \circ T^{-n}
  =\\
  (-1)^k \int_{T^n(W_j)} \prod_{i=k+1}^\smp
  w_{\pi(i)}^u (h_n\circ T^{-n}) \cdot \prod_{i=k}^1 w_{\pi(i)}^s \vf \cdot
  \rho_j
  \circ T^{-n} + O(\|h\|_{p-1,q+1}).
  \end{multline}
Let $\bar w_i^u (x)= DT^n (x)^{-1} w_i^u(T^n x)$. This is a vector
field on a neighborhood of $W_j$. Changing variables, the last
integral in \eqref{eq:kok} is equal to
  \begin{equation}
  \label{eq:koko}
  \int_{W_j}
  \prod_{i=k+1}^\smp
  \bar w_{\pi(i)}^u h_n \cdot \left(\prod_{i=k}^1 w_{\pi(i)}^s
  \vf\right)\circ T^n \cdot \rho_j\cdot J_WT^n,
  \end{equation}
where $J_WT^n$ is the jacobian of $T^n :W_j \to W$, as in the proof of
\eqref{eq1}.

We use the standard coordinate chart (of dimension $d_s$) on $W$,
and write $w_{\pi(i)}^s = \sum_{l=1}^{d_s} g_{\pi(i),l}
\partial_l$ where $g_{\pi(i),l}$ is $\Co^{r+1}$.
Differentiating one of the functions $g_{\pi(i),l}$ yields another term
bounded by $C_n\|h\|_{p-1,q+1}$. Consequently, up to
$O(\|h\|_{p-1,q+1})$, \eqref{eq:koko} is equal to the sum, for
$l_{1}, \dots,l_k
  \in \{1,\dots,d_s\}^{k}$, of
  \begin{equation*}
  \int_{W_j}
  \prod_{i=k+1}^\smp
  \bar w_{\pi(i)}^u  h_n \cdot \left(\prod_{i=1}^k \partial_{l_i}
  \vf\right)\circ T^n \cdot \rho_j \cdot J_WT^n
  \cdot \prod_{i=1}^k g_{\pi(i),l_i} \circ T^n.
  \end{equation*}
Fix parameters $l_1,\dots,l_k$. Let $F= \rho_j \cdot \prod_{i=1}^k
g_{\pi(i),l_i} \circ T^n$. By Lemma \ref{distortion},
$|F|_{\Co^{p+q}(W_j)} \leq C$. We want to estimate
  \begin{equation}
  \label{eq:kokok}
  \int_{W_j}
  \prod_{i=k+1}^\smp \bar w_{\pi(i)}^u  h_n \cdot \left(\prod_{i=1}^k \partial_{l_i}
  \vf\right)\circ T^n \cdot
  J_WT^n \cdot F.
  \end{equation}

Assume first that $\smp=p$ and $k=0$. The function $\vf \circ T^n$
satisfies $|\vf \circ T^n|_{\Co^{p+q}(W_j)} \leq C$. Since
$|F|_{\Co^{p+q}(W_j)} \leq C$, and the vector fields
$w_{\pi(i)}^u$ have a $\Co^{p+q}$ norm bounded by $ \lambda^{-n}$
by Lemma \ref{distortion}, \eqref{eq:def_sharp}
(applied with $f=|\det DT^n|$) implies  that \eqref{eq:kokok} is
bounded by $C \lambda^{-pn} \|h\|_{p,q} \bigl| |\det
DT^n|^{-1}\bigr|_{\Co^{p+q}(W_j)} |J_W T^n
\bigr|_{\Co^{p+q}(W_j)}$.

In the other cases, $\smp-k<p$.  It will be useful to smoothen the
test function. For $\ve\leq \delta$ and $\bar \vf \in
\Co^{q+\smp-k}_0(W,\R)$, let $\A_\ve \bar \vf\in
\Co^{q+1+\smp-k}_0(\tilde{W}, \R)$ be obtained by convolving $\bar
\vf$ with a $\Co^\infty$ mollifier whose support is of size $\ve$.

\begin{lem}\label{lem:ave}
For each $\bar \vf\in\Co^{q+\smp-k}$,
  \[
  \begin{split}
  &|\A_\ve\bar \vf|_{\Co^{q+\smp-k}}\leq C |\bar
  \vf|_{\Co^{q+\smp-k}};\quad\quad |\A_\ve\bar
  \vf|_{\Co^{q+1+\smp-k}}\leq
  C\ve^{-1}|\bar \vf|_{\Co^{q+\smp-k}};\\
  &|\A_\ve\bar \vf-\bar \vf|_{\Co^{q+\smp-k-1}}\leq C\ve |\bar
  \vf|_{\Co^{q+\smp-k}}.
  \end{split}
  \]
\end{lem}
The proof of the above lemma is standard and is left to the reader.

We apply this lemma to $\bar \vf= \prod_{i=1}^k \partial_{l_i}
  \vf$, with $\ve=\nu^{(q+\smp-k)n}$. Then
  \begin{equation*}
  | (\A_\ve \bar \vf - \bar \vf) \circ T^n|_{\Co^{q+\smp-k}(W_j)} \leq C
  \nu^{(q+\smp-k)n}.
  \end{equation*}
Hence, by \eqref{eq:def_sharp},
  \begin{multline*}
  \left| \int_{W_j}
  \prod_{i=k+1}^\smp \bar w_{\pi(i)}^u  h_n \cdot \left(\prod_{i=1}^k \partial_{l_i}
  \vf- \A_\ve \prod_{i=1}^k \partial_{l_i}
  \vf \right)\circ T^n \cdot
  J_WT^n \cdot F\right|
  \\
  \leq C \nu^{(q+\smp-k)n} \lambda^{-(\smp-k)n}
  \bigl| |\det DT^{n}|^{-1}\bigr|_{\Co^{q+\smp-k}(W_j)} \cdot
  |J_WT^n|_{\Co^{q+\smp-k}(W_j)} \|h\|_{p,q}.
  \end{multline*}
The worst bound is obtained when $k=\smp$, in which case
$\nu^{(q+\smp-k)n} \lambda^{-(\smp-k)n}=\nu^q$. Moreover, since
$\smp-k<p$ and  $\A_\ve \bar\vf$ is smoother than $\bar \vf$,
  \begin{equation*}
  \int_{W_j}
  \prod_{i=k+1}^\smp \partial_{l_i}  h_n \cdot
  \left(\A_\ve \prod_{i=1}^k \partial_{l_i}
  \vf \right)\circ T^n \cdot J_WT^n \cdot F
  =O(\|h\|_{p-1,q+1}).
  \end{equation*}
To sum up, we have proved that
  \begin{multline*}
  \left| \int_W v_1 \dots v_\smp(\Lp^n  h) \cdot \vf \right|
  \leq O( \|h\|_{p-1,q+1} )
  \\
  +\left(\sum_j \bigl|
  |\det DT^n|^{-1}\bigr|_{\Co^{r}(W_j)} \cdot |J_W T^n|_{\Co^{r}(W_j)} \right)
  \begin{cases}
  C (\nu^q)^n \norm{h}_{p,q} &\text{if }\smp<p.
  \\
  C \max(\lambda^{-p},\nu^{q})^n \norm{h}_{p,q}  & \text{if }\smp=p.
  \end{cases}
  \end{multline*}
By Lemma \ref{lem:sumdet}, the sum $\sum_j \bigl| |\det
DT^n|^{-1}\bigr|_{\Co^{r}(W_j)} \cdot |J_W T^n|_{\Co^{r}(W_j)}$ is
bounded independently of $n$. This concludes the proof of Lemma
\ref{main_lemma_LY}.
\end{proof}

We now prove \eqref{eq2} by induction over $p$. The case $p=0$ is
given by \eqref{eq1}.

Lemma \ref{main_lemma_LY} implies the inequality
  \begin{equation}
  \label{eq:almost_eq2}
  \| \Lp^n h\|_{p,q} \leq C \max(\lambda^{-p},\nu^{q})^n
  \|h\|_{p,q}+C_n \|h\|_{p-1,q+1}.
  \end{equation}
To prove \eqref{eq:almost_eq2}, we have only
used the fact that $\nu$ is 
greater than the minimal contraction of $T$ in the stable direction, and
$\lambda$ is less than the minimal expansion in the unstable
direction. Let $\lambda'>\lambda$ and $\nu'<\nu$ satisfy the same
conditions, we get in the same way
  \begin{equation}
  \label{eq:almost_final}
  \| \Lp^n h\|_{p,q} \leq C' \max({\lambda'}^{-p},{\nu'}^{q})^n
  \|h\|_{p,q}+C'_n \|h\|_{p-1,q+1}.
  \end{equation}
Finally, choose $n_0$ such that $C'
\max({\lambda'}^{-p},{\nu'}^{q})^{n_0} \leq
\max(\lambda^{-p},\nu^{q})^{n_0}$. Iterating \eqref{eq:almost_final} for
$n=n_0$ (and remembering that $\| \Lp^m h\|_{p-1,q+1} \leq C
\|h\|_{p-1,q+1}$ by the inductive assumption), we obtain
\eqref{eq2}.

\section{General perturbation results}\label{sec:pert}

It is obvious from the previous discussion that all the results
discussed so far -- and in particular Lemmas \ref{lem:LY} and
\ref{lem:exists_partition} -- hold not only for the map $T$, but
also for any map in a $\Co^{r+1}$ open neighborhood $U$ of $T$, or
for any composition of such maps. We will consider perturbations
of $T$ as described in Section \ref{sec:space}, given by a
probability measure $\mu$ on a space $\Omega$ and functions
$g(\omega,\cdot) \in \Co^{p+q}(X,\R_+)$, and we will assume that
all the random diffeomorphisms $T_\omega$ we consider belong to
the above set $U$. In this section, we will prove Theorem
\ref{thm:spectral-stability}.

\begin{lem}\label{lem:mapdist}
For any map $\tilde{T}\in U$, $p+q<r$, holds
\[
\|\Lp_T h-\Lp_{\tilde{T}}h\|_{p-1,q+1}\leq C d_{\Co^{r+1}}(T,
\tilde{T})\|h\|_{p,q}.
\]
\end{lem}
\begin{proof}
Let $\tilde{W}$ be a full admissible leaf, given by an admissible
graph $\chi\in \Xi_i$ defined on a ball $\ovB(x,A\delta)$. We will
use Lemma \ref{lem:exists_partition} with $\gamma=\gamma_0>1$: there
exists a finite number of $\gamma$-admissible graphs
$\chi_1,\dots,\chi_\ell$, such that $\chi_j$ is defined on a ball
$\ovB(x_j, A\gamma\delta) \subset (-\frac{2r_{i(j)}}{3\gamma},
\frac{2r_{i(j)}}{3\gamma})^{d_s}$ for some index $i(j)$, and such
that the corresponding $\gamma$-admissible leaves cover $T^{-1}(W)$.
Write $\rho_j$ for the corresponding partition of unity.

Take $\tilde{T} \in U$. The projection on the first $d_s$
coordinates of
  \begin{equation*}
  \psi_{i(j)}^{-1}\circ \tilde{T}^{-1}\circ T\circ \psi_{i(j)} \circ
  (\Id, \chi_j)(\ovB(x_j,\gamma A \delta))
  \end{equation*}
contains the ball $\ovB(x_j, A\delta)$ if $U$ is small enough.
Hence, it is possible to define a graph $\tilde{\chi}_j$ on
$\ovB(x_j, A\delta)$ whose image is contained in
$\psi_{i(j)}^{-1}(\tilde{T}^{-1}(\tilde{W}))$. Moreover,
$\tilde{T}^{-1}(W)$ is covered by the restrictions of these graphs
to the balls $\ovB(x_j,\delta)$ if $U$ is small enough. Finally,
$\left| \chi_j
-\tilde{\chi}_j\right|_{\Co^{p+q}(\ovB(x_j,A\delta))} \leq C
d_{\Co^{r}}(T,\tilde{T})$.

Let  $|\alpha|\leq p-1$, $\vf\in
\Co^{q+1+|\alpha|}_0(B(x,\delta),\R)$, and set $\tilde
h_j:=h\circ\psi_{i(j)}$. Then
  \begin{equation}
  \label{eq:LT_1}
  \int_{B(x,\delta)}
  \partial^\alpha ((\Lp_{T}h)\circ \psi_i)(\Id,\chi) \cdot \vf
  =\sum_{|\beta|\leq |\alpha|}\sum_{j=1}^\ell
   \int_{B(x_j,\delta)} \partial^\beta \tilde h_j(\Id,\chi_j)\cdot
   F_{\alpha,\beta,T,j}\cdot \rho_j
  \end{equation}
for some functions $F_{\alpha,\beta,T,j}$ bounded in $
\Co^{q+1+|\beta|}$. The same equation holds for
$\Lp_{\tilde{T}}h$, with $\chi_j$ replaced by $\tilde{\chi}_j$ and
$F_{\alpha,\beta,T,j}$ replaced by a function
$F_{\alpha,\beta,\tilde{T},j}$ satisfying 
$|F_{\alpha,\beta,T,j}-
F_{\alpha,\beta,\tilde{T},j}|_{\Co^{q+|\beta|}} \leq C
d_{\Co^{r+1}}(T,\tilde{T})$.

\begin{comment}
More precisely, $\partial^l(h\circ T^{-1})=\sum_{m\leq l} (\partial^m
h)\circ T^{-1} \cdot \Co^{r+m-l}$ if $T\in \Co^{r+1}$ (proof by
induction on $l\geq 1$). Hence, if $k<p$,
  \begin{align*}
  \int \partial^k(h\circ T^{-1} \cdot \det DT^{-1}) \Co^{q+1+k}&
  =\sum_{l\leq k} \int \partial^l (h\circ T^{-1}) \Co^{r-k+l}
  \Co^{q+1+k}
  \\&
  =\sum_{m\leq l\leq k} \int (\partial^m h) \circ T^{-1}
  \Co^{r+m-l}\Co^{r-k+l}  \Co^{q+1+k}.
  \end{align*}
Using the inequalities $m\leq l\leq k\leq p$ and $p+q<r$, we
easily check that the test function, i.e. $F$, belongs to
$\Co^{q+1+m}$. To estimate $|F-\tilde{F}|$ using
$d_{\Co^r}(T,\tilde{T})$, we can not use anymore any bound on $T$
in $\Co^{r+1}$, only in $\Co^r$. Hence, we loose one exponent, and
we can only prove that $F$ and $\tilde{F}$ are close in
$\Co^{q+m}$.
\end{comment}

For $1\leq j\leq \ell$ and $|\beta|\leq |\alpha|$, we have
  \begin{multline}
  \label{eq:LT_12}
  \left| \int_{B(x_j,\delta)} \partial^\beta \tilde h_j(\Id,\chi_j)
  (F_{\alpha,\beta,T,j}-F_{\alpha,\beta,\tilde{T},j}) \rho_j
  \right|\\
  \leq C \norm{h}_{p,q}
  |F_{\alpha,\beta,T,j}-F_{\alpha,\beta,\tilde{T},j}|_{\Co^{q+|\beta|}}
  \leq C \norm{h}_{p,q} d_{\Co^{r+1}}(T,\tilde{T})
  \end{multline}
and
  \begin{multline}
  \label{eq:LT_13}
  \left| \int_{B(x_j,\delta)} \partial^\beta \tilde h_j(\Id,\chi_j)
  F_{\alpha,\beta,\tilde{T},j} \rho_j -
  \int_{B(x_j,\delta)} \partial^\beta \tilde h_j(\Id,\tilde{\chi}_j)
  F_{\alpha,\beta,\tilde{T},j} \rho_j \right|
  \\
  = \left| \int_{t=0}^1
  \int_{B(x_j,\delta)} D( \partial^\beta \tilde h_j)(\Id,\tilde{\chi}_j+ t(
  \chi_j-\tilde{\chi}_j))\cdot (0, \chi_j-\tilde{\chi}_j)
  F_{\alpha,\beta,\tilde{T},j} \rho_j \right|.
  \end{multline}
When $t$ is fixed, each integral is an integral along an admissible
graph, whence it is at most
  \begin{equation*}
  C \|h\|_{p,q}
  |\chi_j-\tilde{\chi}_j|_{\Co^{q+|\beta|+1}}
  | F_{\alpha,\beta,\tilde{T},j}|_{\Co^{q+|\beta|+1}}
  \leq C \| h\|_{p,q}
  d_{\Co^{r+1}}(T,\tilde{T}).
  \end{equation*}
Integrating over $t$, we get $\eqref{eq:LT_13} \leq C \| h\|_{p,q}
d_{\Co^{r+1}}(T,\tilde{T})$.
Combining this inequality with
Equations \eqref{eq:LT_12} and \eqref{eq:LT_1} yields the
conclusion of the lemma.
\end{proof}
This lemma readily implies that, for any operator $\Lp_{\mu,g}$
satisfying the previous assumptions,
  \begin{equation}
  \label{eq:kl_1}
  \|\Lp_{\mu,g}h-\Lp_{T}h\|_{p-1,q+1}\leq
  C\Delta(\mu,g)\|h\|_{p,q},
  \end{equation}
where $\Delta(\mu,g)$ is defined in \eqref{eq:defDeltamug}.

When $g(\omega,x)=1$, Lemma \ref{lem:LY} applied to compositions
of operators of the form $\Lp_{T_\omega}$ immediately implies that
  \begin{equation}
  \label{eq:determ}
  \|\Lp^n_{\mu,g}h\|_{p,q}\leq C\max(\lambda^{-p},\nu^{q})^n
  \|h\|_{p,q}+C\|h\|_{p-1,q+1},
  \end{equation}
which is sufficient to obtain spectral stability, by \cite{KL}. In
particular this suffices to prove Theorem \ref{thm:spectral-stability}
for deterministic perturbations.

However, in the general case, further arguments are required to obtain
a uniform Lasota-Yorke type inequality:
\begin{lem}\label{lem:estrandom}
For any $M>1$ and any perturbation $(\mu,g)$ of $T$ as above,
there exists a constant $C=C\bigl(M,\int |g(\omega,\cdot)|_{
\Co^{p+q}(X,\R)} \dd\mu(\omega)\bigr)$ such that, for any $n\in
\N$,
  \begin{equation}
  \label{eq:kl_2}
  \|\Lp^n_{\mu,g}h\|_{p,q}\leq C
   \max(\lambda^{-p},\nu^{q})^n \|h\|_{p,q}+C M^n
  \|h\|_{p-1,q+1}.
  \end{equation}
\end{lem}
\begin{proof}
We will prove that
  \begin{equation*}
  \| \Lp_{\mu,g}^n h\|_{0,q} \leq C M^n
  \| h\|_{0,q},
  \end{equation*}
by adapting the proof of equation \eqref{eq1}. The proof of
\eqref{eq:kl_2} in the general case is similar, using the same
ideas to extend the proof of \eqref{eq2}. The only problem comes
from the functions $g(\omega_i,x)$, and a distortion argument will
show that their contribution is small. Let $c=\int
|g(\omega,\cdot)|_{\Co^{q}}\dd\mu(\omega)$. Fix parameters
$\oo_n:=(\omega_1,\dots,\omega_n) \in \Omega^n$. Fix also $\ve>0$.
Write $\tilde{g}_i (x)= g(\omega_i,x)+ \ve \frac{ |g(\omega_i,
\cdot)|_{\Co^{q}}}{c}$.

We will write $T_{\oo_i}= T_{\omega_i} \circ \dots \circ
T_{\omega_1}$. Let $W$ be an admissible leaf, $W_1,\ldots, W_\ell$
a covering of $T_{\oo_n}^{-1}W$ by admissible leaves and
$\rho_1,\dots,\rho_\ell$ a corresponding partition of unity, as in
the proof of \eqref{eq1}. Let also $\vf$ be a $ \Co^q$ test
function. Then
  \begin{multline*}
  \int_W \left(\prod_{i=1}^n g(\omega_i, T_{\omega_i}^{-1} \circ\cdots\circ
  T_{\omega_n}^{-1} x)  \right) \Lp_{T_{\oo_n}} h(x) \vf(x)
  \\
  = \sum_{j=1}^\ell \int_{W_j} |\det DT_{\oo_n}|^{-1} h(x) \left( \prod_{i=1}^n
  g(\omega_i, T_{\oo_{i-1}} x) \right) \vf \circ T_{\oo_n}(x)
  J_WT_{\oo_n}(x) \rho_j(x).
  \end{multline*}
Since $W_j$ is admissible, the last integral can be estimated
using the $\Co^q$ norm of $\prod_{i=1}^n
  g(\omega_i, T_{\oo_{i-1}} x)$. Since $|g(\omega_i,
\cdot)|_{\Co^q} \leq \frac{c}{\ve} \tilde{g}_i(x)$ by definition
of $\tilde{g}_i$, Lemma \ref{lem:distortion} shows
that this norm is bounded by $C \exp( C \frac{c}{\ve})
\prod_{i=1}^n
  \tilde{g}_i (T_{\oo_{i-1}} x)$ for any $x\in W_j$.
Combining this estimate with the distortion arguments of the proof of
\eqref{eq1}, we obtain
    \begin{multline*}
  \left|\int_W \left(\prod_{i=1}^n g(\omega_i, T_{\omega_i}^{-1} \cdots
  T_{\omega_n}^{-1} x)  \right) \Lp_{T_{\omega_n}}\cdots
  \Lp_{T_{\omega_1}} h(x) \vf(x) \right|
  \\
  \leq C \exp\left(C\frac{c}{\ve}\right)
  \norm{h}_{0,q}
  \int_{\tilde W} |\det DT_{\oo_n}|^{-1} \left(\prod_{i=1}^n \tilde{g}_i(
  T_{\omega_i}^{ -1} \cdots
  T_{\omega_n}^{-1} x)  \right).
  \end{multline*}
To estimate this last integral, consider the thickening $Z=
\bigcup_{x\in \tilde{W}} W_\rho^u(x)$, where $W_\rho^u(x)$ is the local
unstable manifold of $T_{\oo_n}$ through $x$. Along this manifold,
the function $\prod_{i=1}^n \tilde{g}_i(
  T_{\omega_i}^{ -1} \cdots
  T_{\omega_n}^{-1} x)$ changes of a multiplicative factor at
most $C \exp\left(C \frac{c}{\ve}\right)$, again by
Lemma \ref{lem:distortion}. Hence,
  \begin{align*}
  \int_{\tilde{W}} |\det DT_{\oo_n}|^{-1} \left(\prod_{i=1}^n \tilde{g}_i(
  T_{\omega_i}^{ -1} \cdots
  T_{\omega_n}^{-1} x)  \right)
  \!\!\!\!\!\!\!\!\!\!\!\!\!\!\!
  \!\!\!\!\!\!\!\!\!\!\!\!\!\!\!
  \!\!\!\!\!\!\!\!\!\!\!\!\!\!\!
  &
  \\&
  \leq C \exp\left(C \frac{c}{\ve}\right)\rho^{-d_u}
  \int_Z  |\det DT_{\oo_n}|^{-1} \left(\prod_{i=1}^n \tilde{g}_i(
  T_{\omega_i}^{ -1} \cdots
  T_{\omega_n}^{-1} x)  \right)
  \\&
  = C \exp\left(C \frac{c}{\ve}\right) \int_{T^{-n}(Z)} \prod_{i=1}^n \tilde{g}_i(
  T_{\omega_{i-1}} \cdots
  T_{\omega_1} x)
  \\&
  \leq C \exp\left(C \frac{c}{\ve}\right) \int_X \prod_{i=1}^n \tilde{g}_i(
  T_{\omega_{i-1}} \cdots
  T_{\omega_1} x).
  \end{align*}
Integrating over all possible values of $\omega$, we finally
obtain
  \begin{multline*}
  \| \Lp_{\mu,g}^n h\|_{0,q}
  \leq  C\|h\|_{0,q} \exp\left(C \frac{c}{\ve}\right)
  \\ \cdot
  \int_X \int_{\Omega^n} \prod_{i=1}^n \left(
  g(\omega_i,T_{\omega_{i-1}} \cdots T_{\omega_1} x)+ \ve
  \frac{ |g(\omega_i, \cdot)|_{ \Co^q}}{c} \right) \dd\mu(\omega_1)
  \dots \dd\mu(\omega_n).
  \end{multline*}
Integrating over $\omega_n$ gives a factor $1+\ve$, since $\int
g(\omega_n, y) \dd\mu(\omega_n)=1$ for any $y$. We can then
proceed to integrate over $\omega_{n-1}, \omega_{n-2},\dots$, and
get
  \begin{equation*}
  \| \Lp_{\mu,g}^n h\|_{0,q}
  \leq  C\|h\|_{0,q}
  \exp\left(C \frac{c}{\ve}\right)
  (1+\ve)^n.
  \qedhere
  \end{equation*}
\end{proof}
The inequalities \eqref{eq:kl_1} and \eqref{eq:kl_2} suffice to
apply \cite{KL}, which implies Theorem
\ref{thm:spectral-stability}.

\section{An abstract perturbation theorem}\label{sec:smoothpert}

Let $\B^0\supset \dots \supset \B^s$ be Banach spaces, $0\in
I\subset \R$ a fixed open interval, and $\{\Lp_t\}_{t\in I}$ a
family of operators acting on each of the above Banach spaces.
Moreover, assume that
  \begin{equation}\label{eq:lypert1}
  \exists M>0, \forall\, t\in I,\quad
  \norm{\Lp_t^n f}_{\B^0} \leq CM^n \norm{f}_{\B^0}
  \end{equation}
and
  \begin{equation}\label{eq:lypert2}
  \exists\, \alpha<M,\;\;\forall\, t\in I,\quad
  \norm{\Lp_t^n f}_{\B^1} \leq C \alpha^n \norm{f}_{\B^1} +
  CM^n \norm{f}_{\B^0}.
  \end{equation}
Assume also that there exist operators $Q_1,\ldots,Q_{s-1}$
satisfying the following properties:
  \begin{equation}\label{eq:Qbound}
  \forall\, j=1,\dots, s-1,\;\; \forall\, i \in [j,s],\quad
  \norm{Q_j}_{\B^i\to \B^{i-j}} \leq C
  \end{equation}
and, setting $\Delta_0(t):=\Lp_t$ and $\Delta_j(t):=\Lp_t -\Lp_0
-\sum_{k=1}^{j-1} t^k Q_k$ for $j\geq 1$,
  \begin{equation}\label{eq:Cksmooth}
  \forall\, t\in I,\; \;
  \forall j=0,\dots s, \;\;
  \forall\, i \in [j,s],
  \quad
  \norm{\Delta_j(t)}_{\B^i \to \B^{i-j}} \leq Ct^j.
  \footnote{In fact, this property is used in the proof only for
  $i=s$, and for $(i,j)=(1,1)$.}
  \end{equation}
These assumptions mean that $\Lp_t$ is a continuous, and even a $\Co^s$
perturbation of $\Lp_0$, but the differentials take their values in
weaker spaces. This setting can be applied to the case of smooth expanding
maps (see \cite{Li3} for the argument limited to the case $s=2$) and to the
transfer operator associated to a perturbation of a smooth Anosov
map as we will see in section \ref{sec:smooth}.

For $\rad >\alpha$ and $\delta>0$, denote by $V_{\delta,\rad}$ the
set of complex numbers $z$ such that $|z|\geq \rad$ and, for all
$1\leq k\leq s$, the distance from $z$ to the spectrum of $\Lp_0$
acting on $\B^k$ is $\geq \delta$.

\begin{thm}
\label{thm_1}
Given a family of operators $\{\Lp_t\}_{t\in I}$
satisfying conditions \eqref{eq:lypert1}, \eqref{eq:lypert2},
\eqref{eq:Qbound} and \eqref{eq:Cksmooth}
and setting
\[
R_s(t):= \sum_{k=0}^{s-1} t^k
  \sum_{\ell_1+\dots+\ell_j=k} (z-\Lp_0)^{-1} Q_{\ell_1} (z-\Lp_0)^{-1} \dots
  (z-\Lp_0)^{-1} Q_{\ell_j} (z-\Lp_0)^{-1},
\]
for all $z\in V_{\delta,\rad}$ and $t$ small enough holds true
  \begin{equation*}
  \norm{ (z-\Lp_t)^{-1} -R_s(t)}_{\B^s \to \B^0} \leq C
  |t|^{s-1+\eta}
  \end{equation*}
where $\eta =\frac{\log(\rad/\alpha)}{\log(M/\alpha)}$.
\end{thm}
Hence, the resolvent $(z-\Lp_t)^{-1}$ depends on $t$ in a
$\Co^{s-1+\eta}$ way at $t=0$, when viewed as an operator from
$\B^s$ to $\B^0$.

Notice that one of the results of \cite{KL} in the present setting reads
 \begin{equation}\label{eq:klest}
  \norm{ (z-\Lp_t)^{-1}- (z-\Lp_0)^{-1}}_{\B^1 \to \B^{0}} \leq C |t|^\eta.
  \end{equation}
Accordingly, one has Theorem \ref{thm_1} in the case $s=1$ where
no assumption is made on the existence of the operators $Q_j$.

\begin{proof}[Proof of Theorem \ref{thm_1}]
Iterating the equation
  \begin{equation*}
  (z-\Lp_t)^{-1}=(z-\Lp_0)^{-1}+ (z-\Lp_t)^{-1}(\Lp_t-\Lp_0)(z-\Lp_0)^{-1},
  \end{equation*}
and setting $A(z,t):=(\Lp_t-\Lp_0)(z-\Lp_0)^{-1}$, it follows
\begin{equation}\label{eq:itera}
\begin{split}
  (z-\Lp_t)^{-1}&=\sum_{j=0}^{s-2}(z-\Lp_0)^{-1}
A(z,t)^j+(z-\Lp_t)^{-1} A(z,t)^{s-1}\\
&=\sum_{j=0}^{s-1}(z-\Lp_0)^{-1}
A(z,t)^j+\left[(z-\Lp_t)^{-1}-(z-\Lp_0)^{-1}\right]A(z,t)^{s-1}.
\end{split}
\end{equation}
Next, for each $j\in\N$ and $a\leq s$, using \eqref{eq:Cksmooth},
we can write
\begin{equation}\label{eq:iteone}
A(z,t)^j=\Delta_{a}(t)(z-\Lp_0)^{-1}A(z,t)^{j-1}+\sum_{\ell=1}^{a-1}t^\ell
Q_\ell(z-\Lp_0)^{-1} A(z,t)^{j-1} .
\end{equation}
For $\epsilon=0$ or $1$, we can then prove by induction the
formula, for all $1\leq m\leq j$
\begin{equation}\label{eq:iterm}
\begin{split}
&A(z,t)^j=\sum_{k=1}^m\sum_{\substack{\ell_1+\cdots+\ell_{k-1}<s-\epsilon\\
\ell_i>0}}
t^{\ell_1+\cdots+\ell_{k-1}}Q_{\ell_1}(z-\Lp_0)^{-1}\cdots\\
&\ \ \ \cdots Q_{\ell_{k-1}}(z-\Lp_0)^{-1}
\Delta_{s-\epsilon-\ell_1-\cdots-\ell_{k-1}}(t)(z-\Lp_0)^{-1} A(z,t)^{j-k}\\
&+\ \ \ \sum_{\substack{\ell_1+\cdots+\ell_{m}<s-\epsilon\\
\ell_i>0}}t^{\ell_1+\cdots+\ell_{m}}Q_{\ell_1}(z-\Lp_0)^{-1}\cdots
Q_{\ell_{m}}(z-\Lp_0)^{-1}A(z,t)^{j-m}
\end{split}
\end{equation}

In fact, for $m=1$ the above formula is just \eqref{eq:iteone} for
$a=s-\epsilon$. Next, suppose \eqref{eq:iterm} true for some $m$,
then by \eqref{eq:iteone} it follows
\[
\begin{split}
Q_{\ell_1}(z-\Lp_0)^{-1}\cdots
Q_{\ell_{m}}&(z-\Lp_0)^{-1}A(z,t)^{j-m}=Q_{\ell_1}(z-\Lp_0)^{-1}\cdots
Q_{\ell_{m}}(z-\Lp_0)^{-1}\\
&\times\bigg[\Delta_{s-\epsilon-\sum_{i=1}^{m}\ell_i}(t)(z-\Lp_0)^{-1}
A(z,t)^{j-m-1}\\
&\ \ \ \ \ + \sum_{\ell_{m+1}=1}^{s-\epsilon-\sum_{i=1}^{m}\ell_i -1}
t^{\ell_{m+1}} Q_{\ell_{m+1}}(z-\Lp_0)^{-1} A(z,t)^{j-m-1}\bigg].
\end{split}
\]
Substituting the above formula in \eqref{eq:iterm} we have the formula
for $m+1$.

We can now easily estimate the terms in which a $\Delta_i$ appears.
In fact, $\|A(z,t)\|_{\B^s}\leq C$, and \eqref{eq:Qbound} and
\eqref{eq:Cksmooth} readily imply that
\[
\|Q_{\ell_1}(z-\Lp_0)^{-1}\cdots
Q_{\ell_k}(z-\Lp_0)^{-1}\Delta_{s-\epsilon-\sum_{j=1}^k\ell_j}(t)
(z-\Lp_0)^{-1}\|_{\B^s\to\B^\epsilon}\leq
C|t|^{s-\epsilon-\sum_{j=1}^k\ell_j}.
\]
The theorem follows then from \eqref{eq:itera}, using
\eqref{eq:iterm} with $\epsilon=0$ and $m=j$ to estimate the terms
$(z-\Lp_0)^{-1}A(z,t)^j$, and \eqref{eq:iterm} with $\epsilon=1$
and $m=s-1$ together with \eqref{eq:klest} to show that
$\norm{\left[(z-\Lp_t)^{-1}-(z-\Lp_0)^{-1}\right]A(z,t)^{s-1}}_{\B^s
\to \B^0} \leq C |t|^{s-1+\eta}$.
\end{proof}

\section{Differentiability results}\label{sec:smooth}

In this section, we prove Theorem \ref{thm:smooth} by Applying Theorem
\ref{thm_1}.
To simplify the exposition, we will abuse notations
and systematically ignore the coordinate charts of the manifold
$X$. As we have carefully discussed in the previous sections, this
does not create any problem.

To start with, let us assume that $\delta_*$ is so small that
$\{T_t\;:\; t\in[-\delta_*,\delta_*]\}$ is contained in the
neighborhood $U$ of $T_0$ in which the estimates of the Lasota-Yorke
inequality hold uniformly.

By Taylor formula we have, for each $f\in\Co^r$ and $s\leq r$,
\begin{equation}
\label{eq:taylor} \Lp_{T_t}f=\sum_{k=0}^{s-1}\frac 1{k!}
\frac{d^k}{dt^k}\Lp_{T_t}f\big|_{t=0}
+\int_0^tdt_1\cdots\int_0^{t_{s-1}}dt_s
\left(\frac{d^s}{dt^s}\Lp_tf\right)(t_s).
\end{equation}
Next, for $1\leq k\leq s-1$,
\begin{equation}
\label{eq:qdef} \frac{d^k}{dt^k}\Lp_{T_t}f(x)\big|_{t=0}
=\sum_{\ell=1}^{k}\sum_{|\alpha|=\ell}J_\alpha(k,t,x)
(\Lp_{T_t}\partial^\alpha f)(x)\big|_{t=0} =:k!Q_kf(x)
\end{equation}
for appropriate functions $J_\alpha(k,t,\cdot)\in\Co^r(X,\R)$.

We are now ready to check the applicability of Theorem
\ref{thm_1}. First of all let us define $\B^i:=\B^{p-1+i,q+s-i}$.
Conditions \eqref{eq:lypert1} and \eqref{eq:lypert2} hold with
$\alpha=\max(\lambda^{-p},\nu^{q+s-1})$ by our choice of
$\delta_*$, and $M=1$ by \eqref{eq:determ}. Moreover, for $1\leq
i\leq s$, the essential spectrum of $\Lp$ acting on $\B^i$ is
contained in $\{|z|\leq \max(\lambda^{-(p-1+i)},\nu^{q+s-i})\}$.
Hence, for $1\leq i \leq s$, $\spectr(\Lp: \B^i\to \B^i)
\cap\{|z|> \max(\lambda^{-p},\nu^q)\}$ is composed of isolated
eigenvalues of finite multiplicity. In particular,
$V_{\delta,\rad}$ is discrete.

From the definition of the norms it follows straightforwardly
that, for each multi-index $\alpha$ with $|\alpha|=j$,
$\partial^\alpha$ is a bounded operator from $\B^{p,q}$ to
$\B^{p-j,q+j}$. From this, Condition \eqref{eq:Qbound} readily
follows. By \eqref{eq:taylor} and \eqref{eq:qdef}, it follows that
$\Delta_k$ is given by the last term in \eqref{eq:taylor}. By the
previous arguments
\[
\left\|\frac{d^k}{dt^k}\Lp_{T_t}(f)\right\|_{p-k,q+k}\leq
C\|f\|_{p,q},
\]
which obviously implies Condition \eqref{eq:Cksmooth}.

\appendix

\section{Distortion estimates}
\label{appendix_distortion}

In this appendix, we prove Lemma~\ref{distortion}. Recall that
$\tilde W_j$ and $\tilde W$ are considered as subsets of
$\R^{d_s}\times \{0\} \subset \R^d$. For $y\in \tilde W$, let
$F(y)=\R^{d_s} \times \{0\}$. This defines a $\Co^\infty$ field of
planes of dimension $d_s$ on $\tilde W$. For $x\in \tilde W_j$,
set also $E(x) = DT^n(x)(\{0\}\times \R^{d_u})$. Let
$\vartheta_\ve$ be a $\Co^{\infty}$ mollifier of size $\ve$ on
$\R^{d_s}$. Define a new field of planes of dimension $d_u$ on
$T^n(W_j)$ by $E_\ve(y)= \int E(T^{-n} y+ z) \vartheta_\ve(z)\dd
z$. It is still uniformly transversal to $F$, and it is
$\Co^{r+1}$ along $T^n(W_j)$, even though $E$ was only $\Co^r$,
thanks to the regularizing effect of $\vartheta_\ve$. Note that a
convolution usually shrinks slightly the domains of definition of
functions. However, at the beginning, our functions are defined on
larger sets $\tilde W$ and $\tilde W_j$. Hence, we can safely
forget about this issue in what follows.

A vector field $v$ with $|v|_{\Co^{r+1}}\leq 1$ can be decomposed
along $T^n(W_j)$ as $v=w^u+w^s$ where $w^u(y)\in E_\ve(y)$ and
$w^s\in F(y)$. We will first estimate the norms of this
decomposition along $W$, and prove that, if $\ve$ is small enough,
  \begin{equation}
  \label{estimate_ws}
  |w^s \circ T^n|_{\Co^{r}(W_j)} \leq C
  \end{equation}
and
  \begin{equation}
  \label{estimate_wu}
  |DT^n(x)^{-1}w^u(T^n x)|_{\Co^{p+q}(W_j)} \leq C
  \lambda^{-n}.
  \end{equation}
Then, the second step of the proof will be to extend this
decomposition to a neighborhood of $T^n(W_j)$ so that the
conclusions of Lemma~\ref{distortion} hold.

We will first estimate the $\Co^r$ norm of $x\mapsto E(x)$ along
$W_j$. If $DT^n(x)=\begin{pmatrix}A^n(x) & B^n(x)\\0 &
D^n(x)\end{pmatrix}$, the projection on $E(x)$ is given by
$\begin{pmatrix}0&U^n\\ 0&\Id\end{pmatrix}$ where
$U^n(x)=B^n(x)D^n(x)^{-1}$. Let $x_0$ be an arbitrary point of
$W_j$, we will work on a small neighborhood of $x_0$. 
For $k=1,\dots,n-1$, let $\chart_k$ be a chart on a
neighborhood of $T^k x_0$ such that $\chart_k(T^k(W_j)) \subset
\R^{d_s} \times \{0\}$ and $D\chart_k(T^k x_0) DT^k
(\{0\}\times \R^{d_u})=\{0\}\times \R^{d_u}$. Since the manifolds $T^k(W_j)$
are uniformly $\Co^{r+1}$ (locally, they are admissible leaves)
and uniformly transversal to $DT^k (\{0\}\times \R^{d_u})$, we can choose such
charts with a uniformly bounded $\Co^{r+1}$ norm. Set also
$\theta_0=\Id$ and $\theta_n=\Id$ (this is coherent with the
previous choices since $W_j$ and $W$ are already assumed to be
subsets of $\R^{d_s} \times \{0\}$).

Let $\tilde{T}_k= \chart_{k} \circ T \circ \chart_{k-1}^{-1}$, and
$\tilde{T}^k=\tilde{T}_{k}\circ \dots \circ \tilde{T}_1$. For
$y\in \R^{d_s}\times\{0\}$, we can write
  \begin{equation*}
  D\tilde{T}_k(y)= \left( \begin{matrix} A_k(y) & B_k(y) \\ 0 & D_k(y)
\end{matrix} \right),
  \end{equation*}
with $| A_k(y)| \leq \nu$, $|D_k(y)^{-1}|\leq \lambda$. For $0\leq
k\leq n-2$, since $B_k(\tilde{T}^k x_0)=0$, we can reduce the
neighborhood of $x_0$ and assume that $|B_k(y)|\leq \nu^k$. For
$x\in \R^{d_s} \times \{0\}$, let
  \begin{equation*}
  D \tilde{T}^k(x) = \left( \begin{matrix} A^k(x) & B^k(x) \\ 0 &
D^k(x) \end{matrix} \right).
  \end{equation*}
If $U^k(x)=B^k(x) D^k(x)^{-1}$, we ca write
  \begin{align*}
  &D^{k+1}(x)^{-1}= D^k(x)^{-1} D_{k+1}(\tilde{T}^k x)^{-1},
  \\
  &U^{k+1}(x)=(A_{k+1}(\tilde{T}^k x) U^k(x) + B_{k+1}(\tilde{T}^k x))
D_{k+1}(\tilde{T}^k x)^{-1}.
  \end{align*}
We have $|D_{k+1}( \tilde{T}^k x)^{-1}| \leq \lambda^{-1}$, and its
derivatives with respect to $x$ are bounded by $C \nu^{-k}$ by uniform
contraction of $\tilde{T}$ along $\R^{d_s}\times\{0\}$. Hence,
  \begin{equation}
  \label{D^n_bounded}
  |D^{k}(x)^{-1}|_{\Co^r} \leq \prod_{\ell=1}^k
  (\lambda^{-1}+C\nu^{-\ell})\leq C \lambda^{-k}.
  \end{equation}
In the same way, since $|B_k(y)|\leq \nu^k$ for $k\leq n-2$ by the
smallness of the neighborhood of $x_0$,
  \begin{equation*}
  | U^{k+1}|_{\Co^r} \leq \bigl( (\nu+ C \nu^k) |U^k|_{\Co^r} + C\nu^k
  \bigr) (\lambda^{-1}+C \nu^k).
  \end{equation*}
This implies $|U^{n-1}|_{\Co^r} \leq C \nu^n$, whence
  \begin{equation}
  \label{U^n_bounded}
  | U^n(x) |_{\Co^r} \leq C.
  \end{equation}

Let $v$ be a $\Co^{r+1}$ vector field on a neighborhood of $T^n(W_j)$.
For $x\in T^n(W_j)$ and $y=T^n(x)$, write $ v(y)=(v_1(y),v_2(y))$
the decomposition of ${v}$ along $\R^{d_s} \times \R^{d_u}$. Then
the decomposition of $v(y)$ in $w^u(y)+w^s(y)$ is given by
  \begin{align*}
  &w^u(y)=
  \left(\left[\int U^n(x+z)\vartheta_\ve(z)\dd z\right]
   v_2(y), v_2(y)\right),
  \\&
  w^s(y)=
  \left( v_1(y)-\left[\int U^n(x+z)\vartheta_\ve(z)\dd z\right]
   v_2(y),0\right).
  \end{align*}
Namely, these vectors satisfy $w^u + w^s=v$,
$w^u $ is
tangent to $E_\ve(y)$ and $ w^s$ is
tangent to $\R^{d_s}\times\{0\}$.

Since the $\Co^r$ norm of $U^n$ is bounded, by
\eqref{U^n_bounded}, this proves \eqref{estimate_ws}. Moreover,
  \begin{multline*}
  DT^n(x)^{-1}  w^u(T^n x)
  \\
  =\left( A_n(x)^{-1} \left[ \int (U^n(x+z)-U^n(x))
  \vartheta_\ve(z) \dd z\right] v_2(T^n x),
  \; D^n(x)^{-1} v_2(T^n x)\right).
  \end{multline*}
Hence, \eqref{D^n_bounded} implies that the $\Co^r$ norm of the
second component is bounded by $C\lambda^{-n}$. On the other hand,
the first component is not necessarily small in the $\Co^r$
topology. However, since $p+q<r$, its $\Co^{p+q}$ norm is bounded
by
  \begin{equation*}
  C | A_n(x)^{-1}|_{\Co^{p+q}(W_j)} \ve^{r-(p+q)},
  \end{equation*}
which can be made arbitrarily small by choosing $\ve$ small
enough. This proves \eqref{estimate_wu}.

We still have to extend $w^s$ and $w^u$ to a neighborhood of
$T^n(W_j)$. Let $\pi :\R^d \to \R^{d_s}$ be the projection on the
first $d_s$ components. A naive idea to extend $w^u$ is to set
  \begin{equation*}
  w_1^u(y)=DT^n( T^{-n}y) DT^n(\pi T^{-n}y)^{-1} w^u( T^n \pi
  T^{-n}y).
  \end{equation*}
In other words, we extend $w^u$ so that the vector field
$DT^n(x)^{-1} w_1^u(T^n x)$ is constant along the vertical planes
$\{\eta\}\times \R^{d_u}$. By \eqref{estimate_wu}, this extension
satisfies
  \begin{equation}
  \label{estimate_w1u}
  | DT^n(x)^{-1} w_1^u(T^n x)|_{\Co^{p+q}(V)} \leq C \lambda^{-n},
  \end{equation}
for some neighborhood $V$ of $W_j$. The vector field $w^u_1$ is
unfortunately only $\Co^r$, which means that we will have to
regularize it.
\begin{lem}
Let $G:\R^d \to \R$ be a $\Co^r$ function whose restriction to
$\R^{d_s}\times \{0\}$ is $\Co^{r+1}$. Then, for every $\ve>0$,
there exists a $\Co^{r+1}$ function $H:\R^d\to \R$ such that
\begin{enumerate}
\item $|H|_{\Co^{r+1}} \leq
C_\ve|G|_{\Co^r}+C|G|_{\Co^{r+1}(\R^{d_s}\times\{0\})}$.
\item The restrictions of $G$ and $H$ to $\R^{d_s}\times \{0\}$
are equal.
\item $|G-H|_{\Co^{p+q}} \leq \ve |G|_{\Co^r}$.
\end{enumerate}
\end{lem}
\begin{proof}
Replacing $G$ and $H$ by $G-G\circ \pi$ and $H-G\circ \pi$, we can
assume without loss of generality that $G=0$ on $\R^{d_s}\times
\{0\}$. Let $H_0$ be obtained by convolving $G$ with a $\Co^\infty$
mollifier $\vartheta_\ve$ of size $\ve$ in $\R^d$. 
Let finally $H=H_0-H_0\circ
\pi$. The first and second conclusions of the lemma are clearly
satisfied by $H$.

The functions $G$ and $H_0$ satisfy $|G-H_0|_{\Co^{p+q}}\leq C
\ve^{r-(p+q)}|G|_{\Co^r}$, which can be made arbitrarily small. To
conclude, we have to prove that the $\Co^{p+q}$ norm of $H_0\circ
\pi$ is arbitrarily small. For $\eta\in \R^{d_s}$, we have
  \begin{equation}
  \label{estimate_H}
  H_0(\eta,0)=\int G(\eta+ \eta',\xi')\vartheta_\ve(\eta',\xi')\dd\eta'\dd\xi'.
  \end{equation}
Since $G=0$ on $\R^{d_s}\times \{0\}$, the $\Co^{p+q}$ norm of the
restriction of $G$ to $\R^{d_s}\times \{\xi\}$ is bounded by $C
\xi^{r-(p+q)} |G|_{\Co^r}$. Together with \eqref{estimate_H},
this implies $|H_0\circ \pi|_{\Co^{p+q}}\leq C
\ve^{r-(p+q)}|G|_{\Co^r}$.
\end{proof}

Applying this lemma to the components of $w^u_1$, we obtain a new
vector field $w^u_2$, which coincides with $w^u$ on
$T^n(W_j)$, belongs to $\Co^{r+1}$, and with $|w^u_1 -
w^u_2|_{\Co^{p+q}}\leq \ve$. Choosing $\ve$ small enough, this
together with \eqref{estimate_w1u} implies
  \begin{equation*}
  |DT^n(x)^{-1} w^u_2(T^nx)|_{\Co^{p+q}(V)}\leq C \lambda^{-n}.
  \end{equation*}
Let finally $w^s=v-w_2^u$, the vector fields $w^s$ and $w^u_2$
satisfy all the conclusions of Lemma~\ref{distortion}.

\end{document}